\documentclass[12pt]{article}
\usepackage{amsmath}
\usepackage{graphicx}
\usepackage{amsfonts}
\usepackage{amssymb}
\newtheorem{theorem}{Theorem}[section]
\newtheorem{lemma}{Lemma}[section]

\newtheorem{definition}{Definition}[section]

\newtheorem{remark}{Remark}[section]

\begin{document}

\title{ ON ONE-PARAMETER FAMILIES OF DIDO RIEMANNIAN PROBLEMS}
\author{Moussa BALDE \thanks{$\,\,\,$INSA de Rouen, Laboratoire de Math\'{e}matiques
AMS-LMI, UPRES-A CNRS 6085, BP 08, Pl. E. Blondel, 76131 Mt. St. Aignan cedex,
France. E-mail : balde@lmi.insa-rouen.fr}}
\date{}
\maketitle
\begin{abstract}
Locally, isoperimetric problems on Riemannian surfaces are sub-Riemannian
problems in dimension $3$. The particular case of Dido problems corresponds to
a class of singular contact sub-Riemannian metrics :\textit{\ }metrics which
have the charateristic vertor field\ as symmetry. We give a classification of
the generic conjugate loci ($i.e$. classification of generic singularities of
the exponential mapping) of a \newline $1$-parameter family of $3$-$d$ contact
sub-Riemannian metrics associated to a $1$-parameter family of Dido Riemannian problems.\medskip
\end{abstract}

\section{Introduction}

\subsection{Sub-Riemannian metrics under consideration}

Let $\mathcal{M}$ be any smooth $3$-dimensional manifold. Let $T\mathcal{M}$
and $T^{\ast}\mathcal{M}$ be respectively its tangent and cotangent bundle and
let $\pi:T^{\ast}\mathcal{M}\rightarrow\mathcal{M}$ be the natural projection.

\noindent A \textit{contact sub-Riemannian metric} on $\mathcal{M}$ is a
couple $\Sigma=(\Delta,g),$ where
\[
\text{ }\Delta=\{\Delta_{q}\}_{q\in\mathcal{M}},\Delta_{q}\subset
T_{q}\mathcal{M},
\]
\newline is a contact structure on $\mathcal{M}$ and
\[
g:\Delta\rightarrow R_{+},
\]
a Riemannian metric on $\Delta.$

\noindent Since $\Delta$ is a contact structure and therefore nonintegrable,
$(\Delta,g)$ defines a distance $d$ on $\mathcal{M}$ \cite{BEL}$.$\smallskip\ 

\noindent A contact sub-Riemannian metric being given on $\mathcal{M},$ there
is a 1-form $\omega$ and a vector field $\nu$ on $\mathcal{M}$, both
determined up to orientation by the conditions:

\begin{enumerate}
\item  The distribution $\Delta$ is the kernel of $\omega$ $(\Delta
=\ker\,\omega),$

\item  the $2$-form $d\omega$ restricted to $\Delta$ is the volume form $V$,

\item  the vector field $Z$ is such that $\omega(Z)=1,\;i(Z)\,d\omega=0.$
\end{enumerate}

\noindent The vector field $Z$ is the \textit{characteristic} \textit{vector
field }of the contact sub-Riemannian structure.\smallskip\ 

\noindent The cotangent bundle has the standard symplectic structure. Consider
the Hamiltonian $\mathcal{H}$ of the metric. It is defined as follows:
associated to $(\Delta,\,g)$, there is a cometric on $T^{\ast}\mathcal{M}:$%
\[
\mathcal{H}(\psi)=\frac{1}{2}\sup_{u\in\Delta\backslash\{0\}}\left(
\frac{\psi^{2}(u)}{\left\|  u\right\|  _{g}^{2}}\right)  ,\text{ where }%
\psi(q)\in T_{q}^{\ast}\mathcal{M}.
\]
On fibers of $T^{\ast}\mathcal{M}$, $\mathcal{H}$ is a positive semi-definite
quadratic form, the kernel of which is the annihilator of $\Delta$.

\smallskip\ 

In the remainder of this paper we will assume that:

\begin{enumerate}
\item  The manifold $\mathcal{M}$ is an open subset of $\mathbb{R}^{3}$,
containing the origin $0,$

\item  the distribution $\Delta$ is specified by an orthonormal frame field
$(F,\,G)$ of the metric $g,$ where $F$ and $G$ are two vector fields defined
on $\mathcal{M},$

\item  for any $q\in\mathcal{M},$ $F(q),$ $G(q),$ and $[F,G](q)$ are linearly
independent (the contact condition).
\end{enumerate}

Hence the Hamiltonian $\mathcal{H}$ is given by:
\[
\mathcal{H}=\frac12\left(  (\psi(F))^{2}+(\psi(G))^{2}\right)  .
\]

\textit{In this paper, we deal with contact sub-Riemannian structures that are
invariant under the action of the one parameter group generated by the
characteristic vector field }$Z$\textit{\ (}$i.e.$\textit{\ sub-Riemannian
structures having a symmetry }$Z$\textit{). The reasons for doing so will
become clear when we show that: to characterize minimizing geodesics for this
class of metrics is a dual sub-Riemannian reformulation of the classical
isoperimetric problem of Dido on Riemannian surfaces.}\smallskip\ 

Let $(M,g)$ be a Riemannian surface. We consider on $M$ the following
\textit{iso-area} problem: First we fix two points $z_{0},\,z_{1}\in M$ and a
curve
\[
\tilde{\iota}:[0,1]\rightarrow M,\tilde{\iota}(0)=z_{0},\tilde{\iota}%
(1)=z_{1}.
\]
We are then faced with the following question:

Which curves
\[
\iota:[0,1]\rightarrow M,
\]
connecting $z_{0}$ and $z_{1},$ such that the area
\[
A=\int_{\Omega}V
\]
of the domain $\Omega$ encircled by $\tilde{\iota}$ and $\iota$ is prescribed,
have minimal Riemannian length?

Denoting such an iso-area problem $(M,g,V),$ it is stated (see for instance
\cite{AG2}) that this problem can be reformulated locally in terms of three
dimensional sub-Riemannian geometry.

More precisely, we can associate a germ of Dido Problems $(M,g,V)_{z_{0}}$
with a germ of an oriented sub-Riemannian structure on $\mathcal{M}$ with a
symmetry $Z$, denoted by:
\[
(\mathcal{M},\Delta,g,Z)_{q_{0}},
\]

where,
\[
\mathcal{M}=M\times\mathbb{R}\text{ }\,\,\,\,\,\text{and}\,\,\,\,\,q_{0}%
=\{z_{0},0\}
\]

Let $V=d\alpha$ and $Z=\dfrac{\partial}{\partial w}$, then
\[
\Delta=\ker(dw+\alpha).
\]
If
\[
\delta:[0,1]\rightarrow\mathcal{M}%
\]
is an admissible curve for $\Delta$ with fixed endpoints
\[
\delta(0)=q_{0}\text{ and }\delta(1)=q_{1}=(z_{1},w_{1}),
\]
hence:

\noindent The sub-Riemannian length of $\delta$ is the Riemannian length of
its projection $\delta^{\ast}$ on $M$ and%

\[
w_{1}=\int_{0}^{1}\alpha(\frac{d\delta^{*}}{d\tau})d\tau.
\]

\noindent Therefore projections of sub-Riemannian length minimizers are
solutions of the iso-area problem $(M,g,V).$

\noindent The sub-Riemannian structures associated with $(M,g,V))$ or simply
$(M,g)$ are called the $IsosR$-structure and are denoted by $IsosR.$\smallskip\ 

\noindent\textit{The main aim of this paper is to classify generic conjugate
loci for }$IsosR$\textit{-metrics.}\smallskip\ 

\noindent For some reasons which will become clear in the next section we will
in fact consider $1$-parameter families of Dido Riemannian problems:
\[
\lambda\in I\subset\mathbb{R}\mapsto(M,g(\lambda))\text{, }I\text{ an
interval.}%
\]

\noindent We will denote the families of $IsosR$-metrics associated with
$(M,g(\lambda))$ by $\mathcal{F}$-$IsosR.$

\subsection{Normal forms for Riemannian metrics on surfaces}

Recall the following result.

\begin{theorem}
Let $(M,g)$ be a Riemannian surface. The metric $g$ has the following normal
form. In normal coordinates with pole $0,$ there is an orthonormal frame
(unique up to rotations and up to the action of $SO(2)$ in $T_{0}M\,$)$,$
$(\mathsf{X},\mathsf{Y})$:%

\begin{equation}
(SNF)\left\{
\begin{tabular}
[c]{l}%
$\mathsf{X}=\frac{\partial}{\partial x}+y\,(\beta\,(y\frac{\partial}{\partial
x}-x\frac{\partial}{\partial y}),$\\
\\
$\mathsf{Y}=\frac{\partial}{\partial y}-x\,(\beta\,(y\frac{\partial}{\partial
x}-x\frac{\partial}{\partial y}).$%
\end{tabular}
\right.  \label{snf}%
\end{equation}

Where $\beta(x,y)$ is a smooth function .\thinspace Moreover%
\[
\beta(0,0)=b_{0}=\frac{1}{6}k_{0}.
\]
The constant $k_{0}$ being the gaussian curvature of $M$ at $0.$
\end{theorem}

\begin{remark}
In the normal coordinates specified above, the $k^{th}$-differentials
$\beta^{k}=D^{k}\beta\,|\,_{0}$ are homogeneous polynomials of degree $k$ in
$x,y$ and define symmetric covariant tensors of degree $k$ on $T_{0}M.$

\noindent These tensors, denoted $\beta_{k},$ do not depend on the orientation
and are invariants of the Riemannian structure.
\end{remark}

\subparagraph{\noindent Decomposition of tensor fields}

Let $S^{k}T^{\ast}M$ denote the bundle of covariant symmetric tensors of
degree $k$ on $TM.$ Due to the action of $SO(2)$ on the typical fiber of $TM,
$ we have the following decomposition of $S^{k}T^{\ast}M$ into isotypic
components \cite{AG2}:%

\[
S^{k}T^{*}M=\oplus_{M}(S^{k}T^{*}M)_{j},
\]

\noindent where $(S^{k}T^{\ast}M)_{j}$ is the component relative to the
$j^{th}$ power of the basic character $e^{i\varphi},$ $(i=\sqrt{-1})$.

\noindent Therefore, if $\beta_{n}\in S^{k}T^{\ast}M$ then:
\[
\beta_{n}=\sum_{l}\beta_{n,j},\text{ where }\beta_{n,j}\in(S^{k}T^{\ast}%
M)_{j}.
\]

\subsection{Lemmas}

Let $M$ be any open subset of $\mathbb{R}^{2}$ . Let $O(M)$ and $J^{n}O(M)$ be
respectively the bundle of orthonormal frames of Riemannian metrics on $M$ and
the vector bundle of $n$-jets of elements of $O(M)$.\newline \noindent We
denote by $\Pi_{0}:O(M)\rightarrow M$ and $\Pi:J^{n}O(M)\rightarrow M$ the
standard projections.

\noindent Let
\[
B^{n}=\oplus\{S^{k}T^{\ast}M\,|\,0\leq k\leq n\},B^{-1}=\{0\}
\]

\noindent For $n\geq0$, let us define the following map:%

\[
\Pi^{n}:J^{n+2}O(M)\rightarrow B^{n-1},
\]

\noindent by its restriction to the fiber $\Pi^{-1}(q_{0})$ of the bundle
$\Pi:J^{n+2}O(M)\rightarrow M$:%

\[
\Pi^{n}(j_{q_{0}}^{n+2}(X),j_{q_{0}}^{n+2}(Y))=\beta^{n-1},
\]

\noindent where $\beta^{n-1}$ is the representative of the tensors
$\beta_{n-1,q_{0}}$ in the unique normal coordinates system at $q_{0}$ such that:%

\[
(X(q_{0}),Y(q_{0}))=(\mathsf{X}(q_{0}),\mathsf{Y}(q_{0})),
\]

\noindent where $(\mathsf{X},\mathsf{Y})$ is the normal form of $(X,Y)$ at
$q_{0}.$

\begin{lemma}
\label{tr1} If we fix $(X,Y)\in O(M),$ then:
\[
\Pi^{n}(X,Y):M\rightarrow B^{n-1},\text{ }\Pi^{n}(X,Y)(q_{0})=\Pi^{n}%
(j_{q_{0}}^{n+2}(X),j_{q_{0}}^{n+2}(Y)).
\]
is a surjective submersion.
\end{lemma}

\noindent The proof of lemma \ref{tr1} is the result of computations not given here.

\begin{definition}
A $1$-parameter family $\mathcal{F}$ of orthonormal frames is given by:
\[
\{(X_{\lambda},Y_{\lambda})\in O(M),\lambda\in I\}.
\]
\end{definition}

\begin{lemma}
\label{tr2} If we fix $(\tilde{X},\tilde{Y})\in\mathcal{F},$%
\[
\tilde{\Pi}^{n}(\tilde{X},\tilde{Y}):M\times I\rightarrow B^{n-1},\text{ }%
\Pi^{n}(\tilde{X},\tilde{Y})(q_{0})=\Pi^{n}(j_{(q_{0},\lambda_{0})}%
^{n+2}(\tilde{X}),j_{(q_{0,\lambda_{0})}}^{n+2}(\tilde{Y}))
\]
is a surjective submersion.
\end{lemma}

\noindent Lemma \ref{tr2} is a consequence of lemma \ref{tr1} , and of the
fact that $J^{n}\mathcal{F}$ submerses into $J^{n}O(M).$

\noindent The algebraic lemmas (lemma \ref{tr1} and lemma\ref{tr2}) are the
main elements in establishing the genericity results here.

\noindent Since we have two lemmas, we have the following two choices: We can
envisage an isolated Dido problem and use lemma \ref{tr1}, or a family of such
problems depending on $1$-parameter by using lemma \ref{tr2}.

\noindent\textit{It is obvious that there are nongeneric situations in the
first case, which become generic in the second. Naturally we will consider
}$1$\textit{-parameter families.}

\subsection{Exponential mapping}

Let us first recall briefly the previous results of \cite{Agr}, \cite{CGK} and
\cite{AEG}, which concern generic contact sub-Riemannian structures in
dimension $3.$

\noindent The following result is stated in \cite{AEG}.

\begin{theorem}
\label{tgnc}Consider a germ at $0$ of contact sub-Riemannian metric
$(\Delta,g)$. There is up to orientation, and up to the action of $SO(2)$ on
$\Delta(0),$ a unique coordinate system (\textit{normal coordinates}) with
respect to which the metric has an orthonormal frame $(\tilde{F},\tilde{G})$
of the following form:
\begin{equation}
(NF)\left\{
\begin{tabular}
[c]{l}%
$\tilde{F}=\frac{\partial}{\partial\tilde{x}}+\tilde{y}\,(\,\tilde{\beta
}(\tilde{y}\frac{\partial}{\partial\tilde{x}}-\tilde{x}\frac{\partial
}{\partial\tilde{y}})+\frac{1}{2}(1+\tilde{\gamma})\frac{\partial}%
{\partial\tilde{w}})$\\
\\
$\tilde{G}=\frac{\partial}{\partial\tilde{y}}-\tilde{x}\,(\tilde{\beta
}\,(\tilde{y}\frac{\partial}{\partial\tilde{x}}-\tilde{x}\frac{\partial
}{\partial\tilde{y}})+\frac{1}{2}(1+\tilde{\gamma})\frac{\partial}%
{\partial\tilde{w}}).$%
\end{tabular}
\right.  \label{fn1}%
\end{equation}
Where $\tilde{\beta}$ and $\tilde{\gamma}$ are smooth functions ; satisfying
the following boundary conditions:
\[
\tilde{\beta}(0,0,\tilde{w})=\tilde{\gamma}(0,0,\tilde{w})=\frac
{\partial\tilde{\gamma}}{\partial\tilde{x}}(0,0,\tilde{w})=\frac
{\partial\tilde{\gamma}}{\partial\tilde{y}}(0,0,\tilde{w})=0.
\]
\end{theorem}

\noindent Indeed, ($NF$) can be obtained even when $\Delta$ is not a contact
structure. In fact the normal coordinates are coordinates w.r.t. local
coordinates charts
\[
(U_{\Gamma},(\tilde{w},\tilde{p},\tilde{q})\circ(\pi\circ\exp H)^{-1}),
\]
where $(\tilde{p},\tilde{q},\tilde{r})$ are dual coordinates in $T^{\ast
}\mathcal{M}$ and $\Gamma$ is a smooth curve transversal to $\Delta.$ In the
contact case, such a curve $\Gamma$ can be taken as the integral curve of
$Z.$\smallskip\ 

In $\Gamma$-normal coordinates, geodesics through $\Gamma$, satisfying
transversality conditions w.r.t. $\Gamma$ are straight lines contained in the
plane $\{\tilde{w}=c\}$, where $c$ is a constant (\cite{AEG}). The set
$C_{\Gamma}^{s}$ of points $x$ of $\mathcal{M},$ that are at a distance $s$ of
$\Gamma,$ form a smooth cylinder for small enough values of $s$.

\subsubsection{Approximation of generic conjugate loci and stability}

In \cite{Agr} and \cite{CGK}, it is proved that for contact sub-Riemannian
metrics, there are in essence two generic situations. In the two cases we have
a representation of the exponential mapping as a \textit{suspension}
\textit{of a mapping} between surfaces.

\noindent The first case happens on a complement of a smooth possibly empty
curve $\mathcal{C}$ of $\mathcal{M}$. In this instance, the $3$-jet\thinspace
$(\exp_{3})\,$is a sufficient jet of the exponential mapping.

\noindent This result follows the fact that the suspension of $(\exp_{3})$ is
a ``Whitney map'': that is a \textit{stable} mapping in the sense of
Thom-Mather (\cite{HW}).

\noindent Hence the $3$-jet of the conjugate loci is sufficient to describe it
and the exponential mapping is determined by the $3$-jet of the metric, in a
neighborhood of its singular locus.\smallskip\ 

\noindent The situation on the curve $\mathcal{C}$ is more complex. The full
conjugate loci $CL_{0}=CL,$ splits into two semi-conjugate loci $CL^{+}$ and
$CL^{-}$ corresponding to $w>0$ or $w<0$ in the normal coordinates.

\noindent The intersection $CL_{w}^{\pm}$ of $CL^{\pm}$ with the planes
$\{w=c\neq0\}$ for small enough values of $\left|  w\right|  $, is a closed
curve with $6$ cusps and self-intersections. Hence the suspension of the
exponential mapping is not a ``Whitney map''.

\noindent To conclude that the exponential mapping approximation is stable, it
is necessary to show that all generic self-intersections are transversal.

\noindent In order to classify the semi-conjugate loci, the authors of
\cite{AEG} define the symbol $S$ of $CL^{\pm}$ in the following manner.

\noindent A symbol $S$ is a sequence of 6 rational numbers, $S=(s_{1}%
,\ldots,s_{6})$ modulo cyclic permutation and reflection. If we follow the
curve $\delta$ starting from a cusp point: $s_{i}$ gives half the number of
self-intersections between the $i^{th}$ and the $(i+1)^{th}$ cusp point.

\begin{theorem}
\label{tc1} $1)$ At generic points of the curve $\mathcal{C}$, the possible
symbols for generic semi-conjugate loci are:%

\[%
\begin{array}
[c]{rrr}%
S_{1}=(2,1,1,2,1,0), &  & S_{2}=(2,1,1,1,1,1),\\
&  & \\
& S_{3}=(0,1,1,1,1,1). &
\end{array}
\]

$2)$ There are two types of isolated points at $\mathcal{C}\frak{.}$

a) At the first type the possible symbols are:
\[%
\begin{array}
[c]{rrr}%
S_{1}^{\ast}=S_{1}=(2,1,1,2,1,0), &  & S_{2}^{\ast}=S_{2}=(2,1,1,1,1,1),\\
&  & \\
& S_{3}^{\ast}=S_{3}=(0,1,1,1,1,1). &
\end{array}
\]

But in that case , the exponential mapping is determined by the higher order
jet of the metric $(7$-jet$)$ than at generic points of $\mathcal{C}$ $(5$-jet$).$

b) At the second type the possible symbols are:
\[%
\begin{array}
[c]{cc}%
S_{4}=(\frac{1}{2},\frac{1}{2},1,0,0,1) & S_{5}=(1,\frac{1}{2},\frac{1}%
{2},1,1,1)\\
& \\
S_{6}=(\frac{3}{2},\frac{1}{2},1,1,0,1) & S_{7}=(2,\frac{1}{2},\frac{1}%
{2},2,0,0).
\end{array}
\]
\end{theorem}

\noindent However, this result, stated also in \cite{AEG}, only gives
classification of possible semi-conjugate loci. The problem of classifying
possible full conjugate loci is therefore not entirely solved.

\subsection{Normal Dido coordinates}

\begin{theorem}
\label{thinf} Given $(\Delta,g)_{0}$ a germ at the origin of an isoperimetric
sub-Riemannian metric. There is an unique coordinates system (normal Dido
coordinate), up to orientation and up to the action of $SO(2)$ on $\Delta
(0),$with respect to which the metric has an orthonormal frame $(\tilde{F}%
_{I},\tilde{G}_{I})$ of the following form:
\begin{equation}
(INF)\left\{
\begin{tabular}
[c]{l}%
$\tilde{F}_{I}=\mathsf{X}+\frac{y}{2}\gamma\frac{\partial}{\partial w}$\\
\\
$\tilde{G}_{I}=\mathsf{Y}-\frac{x}{2}\gamma\frac{\partial}{\partial w}.$%
\end{tabular}
\right.  \label{inf}%
\end{equation}
Where : $(\mathsf{X,Y})$ is the normal form of the Riemannian metric on $M$
and
\begin{equation}
\gamma=(1+(x^{2}+y^{2})\beta)\,\smallint_{0}^{1}\frac{2\,t\,dt}{1+t^{2}%
(x^{2}+y^{2})\,\beta(t\,x,t\,y)}. \label{forg}%
\end{equation}
\end{theorem}

\noindent This theorem is not proved here. Observe that in comparison with
theorem \ref{tgnc}, the only work we have to do for this, is to slightly
modify the $\Gamma$-normal coordinates. This done, checking formula \ref{forg}
is very easy.

\noindent In accordance with \cite{CGK} the $\mathcal{F}$-$IsosR$-structures
are nongeneric. However this is fairly interesting that as it is in the case
of generic contact sub-Riemannian metrics, our work follows almost the same
steps and we prove similar results.

\subsection{Statement of our main results and outline of the paper}

Our work is organized as follows. In the second section, we give the notations
for the main invariants. We end this section by summarizing the basic
properties of the relevant jet of the conjugate locus computed in appendix
\ref{appA2}.

\noindent In section 3, we state our genericity results. We highlight two
principal invariants denoted by $r_{2}$ and $r_{3}.$ More precisely we prove
:\smallskip\ 

\noindent\textit{For a generic element }$\Sigma$\textit{\ of }$\mathcal{F}%
$\textit{\ (for the Whitney topology), the set of points of }$M\times I$
\textit{on which the invariant }$r_{2}$\textit{\ vanishes is a smooth possibly
empty curve }$C$\textit{\ }$.$

\noindent\textit{On the curve }$C$\textit{\ the invariant }$r_{3}%
$\textit{\ does not vanish}.\smallskip\ 

\noindent Using higher order invariants, theorems \ref{th1} and \ref{th2} show
that the curve $C$ is partitioned into two subsets: a discrete subset
(isolated points of $C$) and its complement (generic point of $C$).

\noindent In section 4, we state our stability results (see theorems
\ref{th4c}, \ref{th6c1},\ref{th6c2} and \ref{th6c3}). We can summarize these
results as follow:

\noindent If $r_{2}$ is nonzero, the exponential mapping is determined by the
$5$-jet of the metric, in a neighborhood of its singular locus. . This is the
purpose of theorem \ref{th4c}.

\noindent On $C,$ the situation is by far more complex and particularly
delicate at isolated points on this curve.

\noindent For the sake of clarity theorem \ref{th6c1} deals only with generic
points on the curve $C$. Thus theorems \ref{th6c2} and \ref{th6c3} are devoted
to isolated points on $C$.\smallskip\ 

\noindent It is clear that for $\mathcal{F}$-$IsosR,$ we can classify the semi
conjugate loci as in theorem \ref{tc1}, but we need higher order
jets.\smallskip\ 

\noindent\textit{What is more, for }$\mathcal{F}$\textit{-}$IsosR,$%
\textit{\ our symbols are the same for }$CL^{+}$\textit{\ and }$CL^{-}%
$.\textit{\ It follows that we can classify the full conjugate loci}%
.\smallskip\ 

\noindent In fact we state the following results.

\begin{theorem}
\label{TC}$1)$ At generic points of the curve $C$, the possible symbols for
generic semi-conjugate loci are:%

\[%
\begin{array}
[c]{rrr}%
S_{1}=(2,1,1,2,1,0), &  & S_{2}=(2,1,1,1,1,1),\\
&  & \\
& S_{3}=(0,1,1,1,1,1). &
\end{array}
\]

$2)$ There are two types of isolated points at $C.$

a) At the first type the possible symbols are:
\[%
\begin{array}
[c]{rrr}%
S_{1}^{\ast}=S_{1}=(2,1,1,2,1,0), &  & S_{2}^{\ast}=S_{2}=(2,1,1,1,1,1),\\
&  & \\
& S_{3}^{\ast}=S_{3}=(0,1,1,1,1,1). &
\end{array}
\]

But in that case , the exponential mapping is determined by higher order jets
of the metric $(9$-jet$)$ than at generic points of $C$ $(7$-jet$).$

b) At the second type the possible symbols are:
\[%
\begin{array}
[c]{cc}%
S_{4}=(\frac{1}{2},\frac{1}{2},1,0,0,1) & S_{5}=(1,\frac{1}{2},\frac{1}%
{2},1,1,1)\\
& \\
S_{6}=(\frac{3}{2},\frac{1}{2},1,1,0,1) & S_{7}=(2,\frac{1}{2},\frac{1}%
{2},2,0,0).
\end{array}
\]
\end{theorem}

\begin{theorem}
For a fixed element of $\mathcal{F}$-$IsosR$ the possible symbols for generic
conjugate loci are:%

\[%
\begin{array}
[c]{rrr}%
S_{1}=(2,1,1,2,1,0), &  & S_{2}=(2,1,1,1,1,1),\\
&  & \\
& S_{3}=(0,1,1,1,1,1). &
\end{array}
\]
\end{theorem}

The section 5 is devoted to technical computations of the exponential mapping
and the conjugate loci.

\section{Conjugate loci}

Let $\Sigma=(\Delta,g)_{0}$ be a germ at the origin of an isoperimetric
sub-Riemannian metric$.$ We can restrict $\mathcal{M}$ in order to obtain a
neighborhood $N$ of $x_{0}=(0,0,0)$, with respect to which $\Sigma$ is under
normal form $(INF)$. We denote by $(p,q,r)$ the dual coordinate in $T^{\ast
}N.$ The set $C_{x_{0}}^{\frac{1}{2}}=\pi_{N}^{-1}(x_{0})\cap\mathcal{H}%
^{-1}(\frac{1}{2})$ is therefor the cylinder $\{p^{2}+q^{2}=1\}.$

\smallskip\ 

\noindent Let $\exp sH:T^{\ast}\mathbb{R}^{3}\rightarrow T^{\ast}%
\mathbb{R}^{3}$, denotes the Hamiltonian flow of $H$ at time $s.$ If $N$ is
small enough, point of the conjugate loci in $N$ appear only along geodesics
$\{\pi\circ\exp sH(\lambda_{x_{0}})\}$ for $\lambda_{x_{0}}$ in a certain
neighborhood of infinity $N^{\ast}\subset C_{x_{0}}^{\frac{1}{2}}\,$\cite{AEG}$.$

\smallskip\ 

\noindent In appendix \ref{appA2} we compute relevant jets of $CL$ in a
suitable coordinates system $(h,\varphi),$ with respect to which:
\[
CL:N^{\ast}\rightarrow N
\]
has the following expansion with respect to $h$:
\[
CL(\varphi,\,h)=(x(\varphi,\,h),\,y(\varphi,\,h),\,h)=(\sum_{i=4}^{9}%
f_{i}(\varphi)\,h^{i}+O(h^{10}),\,h),
\]
where $h=\sqrt{\frac{\epsilon w}{\pi}},$ $\epsilon=\pm1\,$ according to the
sign of $w,$ and $(x,y,w)$ the normal coordinates.
\[
CL=CL^{+}\cup CL^{-}%
\]

\smallskip\ 

\noindent As we will see, for small enough of $\left|  w\right|  $, the
intersection of $CL^{\pm}$ with planes $\{w=c\neq0\}$ will be a close curve.
In the more degenerate generic cases this curve will have self-intersections
denoted by $Isoself.$

\subsection{Notation}

The following tensors have the leading part in our study.

\noindent Let $\beta_{1}=\beta_{1,1},$ $\beta_{2}=\beta_{2,2}+\beta_{2,0},$
and $\beta_{3}=\beta_{3,3}+\beta_{3,1}.$ Thus we have in our normal coordinates%

\[%
\begin{tabular}
[c]{l}%
$\beta_{1,1}=R_{e}(r_{1}(dx+idy)).$\\
\\
$\beta_{2,2}=R_{e}(r_{2}(dx+idy)^{2}),$\thinspace$\,\,\beta_{2,0}=\tau
_{0}\,(dx^{2}+dy^{2}).$\\
\\
$\beta_{3,3}=R_{e}(r_{3}(dx+idy)^{3}),$ $\beta_{3,1}=R_{e}(\bar{v}%
\,(dx+idy))(dx^{2}+dy^{2}).$%
\end{tabular}
\]

\noindent Where :%

\[%
\begin{tabular}
[c]{l}%
$r_{1}=\left|  r_{1}\right|  (\cos\theta_{1},\sin\theta_{1})=\left|
r_{1}\right|  (\cos\theta_{1}+i\sin\theta_{1})$\\
\\
$r_{2}=\left|  r_{2}\right|  (\cos\theta_{2},\sin\theta_{2})=\left|
r_{2}\right|  (\cos\theta_{2}+i\sin\theta_{2})$\\
\\
$r_{3}=\left|  r_{3}\right|  (\sin\theta_{3},-\cos\theta_{3})=\left|
r_{3}\right|  (\sin\theta_{3}-i\cos\theta_{3})$\\
\\
$v=-v_{1}+i\,v_{2}.$%
\end{tabular}
\]

\noindent The remaining tensors are given in appendix \ref{appA2}.

\subsection{Basic properties of $f_{i}$ $(B.P.)\,\label{bf}$}

\noindent From now on, we summarize basic properties of the map $f_{i}$ when
$\epsilon=1.$ For details see appendix \ref{appA2}

\begin{itemize}
\item [ P$_{1}$)]$r_{2}\neq0.$

\begin{enumerate}
\item $\,f_{4}(\varphi)=3\pi\left|  r_{1}\right|  (-\sin\theta,\cos\theta),$

\item $f_{4}$ is independant of $\varphi,$

\item $f_{5}(\varphi)=5\pi\,\left|  r_{2}\right|  \,(3\cos(\varphi-\theta
_{2})+\cos(3\varphi-\theta_{2}),3\sin(\varphi-\theta_{2})-\sin(3\varphi
-\theta_{2})).$
\end{enumerate}

\item[ P$_{2}$)] $r_{2}=0.$

\begin{enumerate}
\item $f_{4}(\varphi)=3\pi\left|  r_{1}\right|  (-\sin\theta_{1},\cos
\theta_{1}),\,\,\,\,\,\,f_{5}(\varphi)=(0,0),$

\item
\[%
\begin{array}
[c]{c}%
f_{6}(\varphi)=\frac{\pi}{2}(-25v_{2}+90\left|  r_{3}\right|  \cos
(2\varphi+\theta_{3})+45\left|  r_{3}\right|  \cos(4\varphi+\theta
_{3})+31b_{0}\left|  r_{1}\right|  \sin\theta_{1},\\
\\
-25v_{1}-90\left|  r_{3}\right|  \sin(2\varphi+\theta_{3})+45\left|
r_{3}\right|  \sin(4\varphi+\theta_{3})-31b_{0}\left|  r_{1}\right|
\cos\theta_{1}),
\end{array}
\]

\item $f_{6}(\varphi+\pi)=f_{6}(\varphi),$

\item $f_{7}(\varphi+\pi)=-f_{7}(\varphi),$

\item
\[
\frac{df_{6}}{d\varphi}\wedge\frac{df_{7}}{d\varphi}=12960\pi^{3}\left|
r_{3}\right|  \sin(3\varphi+\theta_{3})\left|  r_{1}^{2}\right|  \sin
2(\varphi-\theta_{1}),
\]

\item
\[
\frac{df6}{d\varphi}\wedge f_{7}=1080\pi^{2}\left|  r_{3}\right|
\sin(3\varphi+\theta_{3})P(\varphi)=\Psi(\varphi),
\]

\item $f_{8}(\varphi+\pi)=f_{8}(\varphi),$

\item $f_{9}(\varphi+\pi)=-f_{9}(\varphi),$

\item
\[
\frac{df_{6}}{d\varphi}\wedge\frac{df_{8}}{d\varphi}=\left|  r_{3}\right|
\sin(3\varphi+\theta_{3})\,(Q_{0}(\varphi)),
\]

\item
\[
\frac{df_{6}}{d\varphi}\wedge f_{9}=\left|  r_{3}\right|  \sin(3\varphi
+\theta_{3})\,(Q_{1}(\varphi)+Q_{2}(\varphi))).
\]
\end{enumerate}
\end{itemize}

\begin{remark}
$P(\varphi)=A\,\cos(2\varphi)+B\,\sin(2\varphi)+C\,\cos(4\varphi
)+D\,\sin(4\varphi),$ where $A,\,B,\,C,\,D$ are linear combinations of the
coefficient of the invariants $\beta_{4}$ and $r_{1}^{2}.$\newline Contrary to
the generic contact case, $P(\varphi)$ is indepent on $\epsilon$ (see appendix
\ref{appA2}).
\end{remark}

\begin{remark}
\label{rp692}\thinspace The $Q_{i}$'s do not depend on the same invariants.
\end{remark}

\noindent In particular only $Q_{2}$ is dpendent on the coefficients of
$\beta_{6}.$

\noindent Let :%
\[
\nu=\dfrac{A-i\,B}{2},\;\mu=\dfrac{C-i\,D}{2}.
\]

\noindent Let
\[%
\begin{tabular}
[c]{l}%
$P(\varphi)=\nu\,e^{2\,i\,\varphi}+\overline{\nu}\,e^{-2\,i\,\varphi}%
+\mu\,e^{4\,i\,\varphi}+\overline{\mu}\,e^{-4\,i\,\varphi}.$\\
\\
and\\
\\
$\widetilde{P}(z)=\mu\,z^{4}+\nu\,z^{3}+\overline{\nu}\,z+\overline{\mu}$\\
\\
$T_{c}(\varphi)=\left|  r_{3}\right|  \sin(3\varphi+\theta_{3})$\\
\\
and\\
\\
$T(z)=r_{3}\,z^{3}+\bar{r}_{3}.$%
\end{tabular}
\]

\noindent In the remainder of this paper we will be concerned with the roots
on the unit circle of these trigonometric polynomials.

\section{Genericity results}

\noindent In the next section we will prove the following results :

\begin{theorem}
\label{th0}

For a generic element of $\mathcal{F}$ (for the Whitney topology), the set on
which $\beta_{2,2}$ vanishes is a smooth curve $C$ of $\mathbb{R}^{3}.$

Along $C,$ $\beta_{3,3}$ and $\mu$ do not vanish$.$
\end{theorem}

\begin{theorem}
\label{th1} On the complement of a discrete subset of $C$ :

i ) $\widetilde{P}(z)$ has either two or four simple roots on the unit circle$,$

ii ) $\widetilde{P}(z)$ and $T(z)$ have no common roots on the unit circle.
\end{theorem}

\begin{theorem}
\label{th2} At isolated point of $C,$ either :

i) $\widetilde{P}(z)$ has a double root on the unit circle, which is not
triple, and which is not a root of $T(z).$

The other roots of $\widetilde{P}(z)$ are simple ; they are on the unit circle
and are not roots of $T(z).$

or,

ii) $\widetilde{P}(z)$ and $T(z)$ have one and only one common simple root on
the unit circle.

The other roots of $\widetilde{P}(z)$ are simple.
\end{theorem}

\subsection{Proof of our genericity theorems}

\subsubsection{The bad set}

\begin{definition}
$B_{0}$ is the subset of $B^{6}$\textrm{\ }such that :

1) $\beta_{2,2}=\beta_{3,3}=0\smallskip$

$B_{1}$ is the subset of $B^{6}$\textrm{\ }such that :

1) $\beta_{2,2}=\mu=0$
\end{definition}

\begin{definition}
$B_{2}$ is the subset of $B^{6}$\textrm{\ }such that :

1) $\beta_{2,2}=0$

2) $\widetilde{P}(z)$ has a triple root on the unit circle.\smallskip

$B_{3}$ is the subset of $B^{6}$ such that :

1) $\beta_{2,2}=0,$

2) A root of $T(z)$ is a double root of $\widetilde{P}(z)$ on the unit
circle\smallskip.

$B_{4}$ is the subset of $B^{6}$\textrm{\ }such that :

1) $\beta_{2,2}=0,$

2) $T(z)$ and $\widetilde{P}(z)$ have two common roots on the unit circle.
\end{definition}

\begin{definition}
$B_{5}$ is the subset of $B^{6}$\textrm{\ }such that

1) $\beta_{2,2}=0,$

$\widetilde{P}(z)$ has a double root on the unit circle\smallskip

$B_{5}$ is the subset of $B^{6}$\textrm{\ }such that

1) $\beta_{2,2}=0,$

$\widetilde{P}(z)$ and $T(z)$ have one and only one common simple root on the
unit circle

($i.e.\operatorname{Re}$s$(r_{3},\mu,\nu)=$Resultant of $\widetilde{P}(z)$ and
$T(z)$ has a simple root).
\end{definition}

\begin{remark}
Since the integral of $\widetilde{P}(z)$ over its period\thinspace is zero,
then, if $\mu\neq0,$ the following facts are easy to check.

1) $\widetilde{P}(z)$ has a root on the unit circle.

2) $\widetilde{P}(z)$ does not have two double roots on the unit circle.

3) If $\widetilde{P}(z)$ has a double root on the unit circle, the two other
roots of $\widetilde{P}(z)$ are on the unit circle.
\end{remark}

\subsubsection{Estimate of the codimension of B$_{k}$'s}

Some appropriate computations on $T(z)$ and $\widetilde{P}(z)$ allow us to
rewrite the $B_{k}$'s as follows:%

\begin{equation}%
\begin{tabular}
[c]{l}%
$B_{0}=\{r_{2}=0\}\cap\{r_{3}=0\}$\\
\\
$B_{1}=\{r_{2}=0\}\cap\{\mu=0\}$\\
\\
$B_{2}=\{r_{2}=0\}\cap\{4\mu^{2}\bar{\nu}+\nu^{3}=0\}$\\
\\
$B_{3}=\{r_{2}=0\}\cap\{27\nu^{3}r_{3}\bar{r}_{3}^{2}+(\bar{\nu}r_{3}-4\mu
\bar{r}_{3})^{3}=0\}$\\
\\
$B_{4}=\{r_{2}=0\}\cap\{\bar{\nu}r_{3}-\mu\,\bar{r}_{3}=0\}$\\
\\
$B_{5}=\{r_{2}=0\}\cap\{27\,\operatorname{Re}^{2}(\mu\bar{\nu}^{2})-(4\left|
\mu\right|  ^{2}-\left|  \nu\right|  ^{2})^{3}=0\}$\\
\\
$B_{6}=\{r_{2}=0\}\cap\{\text{Res(}\mu,\nu,r_{3}\text{)}=0\}$%
\end{tabular}
\label{codi1}%
\end{equation}

\subsubsection{Proof of the genericity}

We consider the $\hat{B}_{k}$'s defined from the $B_{k}$'s by :%

\[
\hat B_{k}=(\tilde\Pi_{\mathcal{N}}^{7})^{-1}(B_{k})
\]

Now let :
\[
\tilde B_{k}=\{(\tilde X,\tilde Y)\,|\,J^{9}(\tilde X,\tilde Y)\pitchfork\hat
B_{k}\}
\]

If $B=\cap_{k}\tilde B_{k}$ ; then $B$ is an open dense set for the Whitney
topology in $\mathcal{F}$.

\noindent Standard arguments from transversality theory (see for instance
Goresky-MacPherson \cite{GM}), allow us to conclude that :%

\[
\forall\,(\tilde X,\tilde Y)\in B,\,(\tilde\Pi_{\mathcal{N}}^{n}(\tilde
X,\tilde Y)^{-1})(B_{k})
\]

are smooth submanifolds or Whitney-stratified set of $M\times I$ of the same
codimension as the $B_{k}$'s.

\noindent Since $B_{k}$'s. are at least codimension $4$ (resp codimension $3)$
for $0\leq k\leq4$ (resp $5\leq k\leq6$), this ends the proofs of theorems
\ref{th0}, \ref{th1}.and \ref{th2}\smallskip.

\begin{remark}
\ If we consider a fixed element $\mathcal{F}_{\lambda_{0}}$ of $\mathcal{F}$,
then: for a generic element of $\mathcal{F}_{\lambda_{0}}$ (for the Whitney
topology), $\beta_{2,2}$ vanishes at isolated points of $\mathbb{R}^{2}.$
\end{remark}

\begin{remark}
\label{RC} at the isolated points of $C$ specified in the part $(i)$ of the
theorem \ref{th2}, we have :%

\[%
\begin{tabular}
[c]{l}%
$Q_{i}\neq0,\,i=0,1,2$\\
\\
$Q_{0}+Q_{2}\neq0$%
\end{tabular}
\]
\end{remark}

Denoting $\hat{B}_{5}$ the subset of $B^{6}\times S^{1}$ (where $S^{1}$ is the
unit circle) defined by the following equations :%

\[
\beta_{2,2}=0,P(\varphi)=0,\frac{dP}{d\varphi}(\varphi)=0,Q_{i}(\varphi
)=0\text{ }(\text{resp }Q_{0}+Q_{2}=0).
\]

Taking into account the remark \ref{rp692}, we show as in \cite{AEG} that the
projection $B_{5}$ of $\hat B_{5}$ on $B^{6}$ is a semi-algebraic subset of
codimension four.

\section{Stability results}

Let $CL$ be the conjugate locus mapping for a generic element of $\mathcal{F}
$-$IsosR$ (for the Whitney topology) and $\mathbb{S}$ a neighborhood of the
singular locus at the source.

\begin{theorem}
\label{th4c} On the complement of the smooth curve $C$ $(\beta_{2,2}\neq0).$

There is a neighborhood $U$ of $\mathbb{S}\cap\{0<h<c\}$ (for small enough
values of $c$), such that :

The $5$-$jet$ of $CL$ : $CL^{5}=f_{4}h^{4}+f_{5}h^{5}$ is a sufficient jet of
$CL$ on $U.$

The restriction of the exponential mapping $\exp|U$ is $5$-determined in $h.$

$\exp|U$ is determined by the $5$-$jet$ of the metric.
\end{theorem}

\noindent The main argument for the proof of the theorem \ref{th4c} is the
fact that the suspension of the exponential mapping is a \ ``Whitney map''

\noindent In fact the intersection $CL_{w}^{\pm}$ between $CL$ and the planes
$w=c$ (for small enough values of $\left|  w\right|  $) is a closed curve
which has fold-points, cusp-points (four) and without self-intersection (see
figures 1 \ and 2).

\noindent Hence the mapping%
\[
f_{4}h^{4}+f_{5}h^{5}+O(h^{6})
\]
is R.L-equivalent to
\[
f_{4}h^{4}+f_{5}h^{5}%
\]

\paragraph{Figures}%

\begin{center}
\includegraphics[
trim=0.000000in 0.000000in 0.002402in 0.002393in,
height=10.0254cm,
width=13.4148cm
]%
{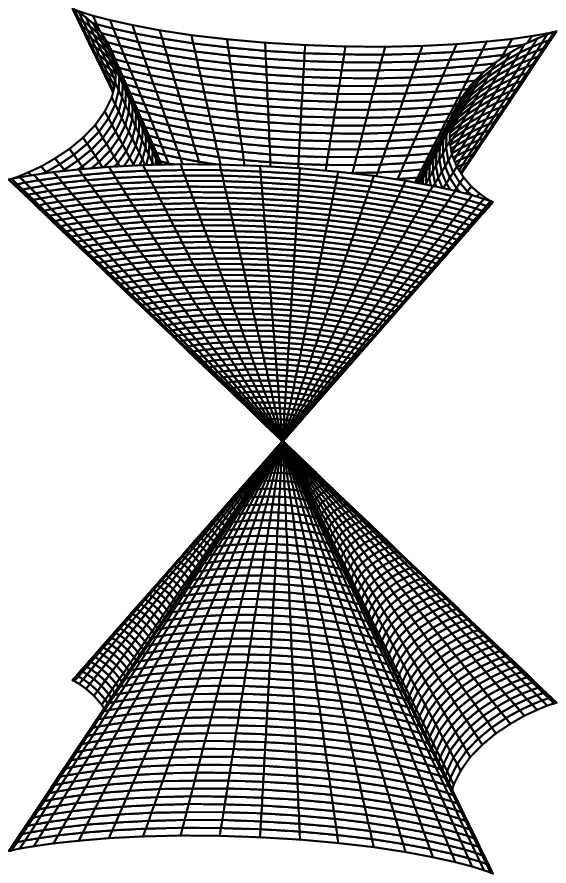}%
\\
$CL^{5}\,\ $when $r_{2}\neq0$.
\end{center}

\begin{center}
\includegraphics[
height=10.0979cm,
width=13.4895cm
]%
{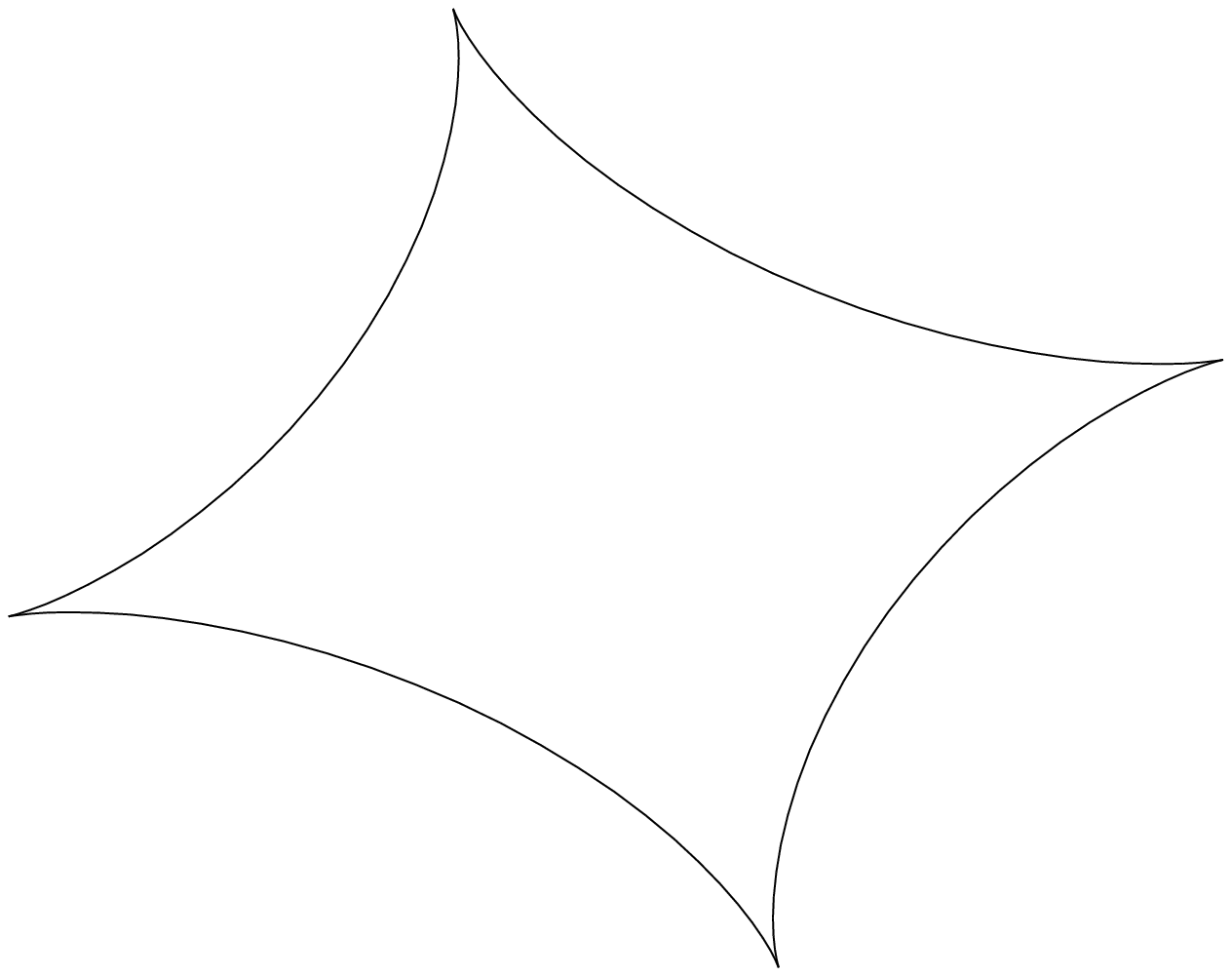}%
\\
$CL_{w}^{\pm}$ \ \ (when $r_{2}\neq0$)
\label{CLP}%
\end{center}
\newpage

\begin{theorem}
\label{th6c1} There is an open dense subset $O$, (complement of a discrete
subset) of $C$, such that :
\[
Isoself=Isoself_{c}+Isoself_{1}%
\]

$-Isoself_{c}$ is the union of three curves.

$-Isoself_{1}$ is the union of two or four curves.

All the self-intersections are transversal and are not dependent on $\epsilon.$

There is a neighborhood $U$ of $\mathbb{S}\cap\{0<h<c\}$ (for small enough
values of $c$), such that :

The $7$-$jet$ of $CL$ : $CL^{7}=f_{4}h^{4}+f_{6}h^{6}+f_{7}h^{7}$ is a
sufficient jet of $CL$ on $U.$

The restriction of the exponential mapping $\exp|U$ is $7$-determined in $h.$

$\exp|U$ is determined by the $9$-$jet$ of the metric.
\end{theorem}

At the isolated points of $C$ we have :

\begin{theorem}
\label{th6c2}

If $P(\varphi)$ has a double root on the circle unity or $P(\varphi)=0$ and
$f_{7}(\varphi)=0$, then :

$-Isoself_{c}$ is the union of three curves.

$-Isoself_{1}$ is the union of two or four curves.

All the self-intersections are transversal and are dependent on $\epsilon.$

There is a neighborhood $U$ of $\mathbb{S}\cap\{0<h<c\}$ (for small enough
values of $c$), such that :

The $9$-$jet$ of $CL$ : $CL^{9}=f_{4}h^{4}+f_{6}h^{6}+f_{7}h^{7}+f_{8}%
h^{8}+f_{9}h^{9}$ is a sufficient jet of $CL$ on $U.$

The restriction of the exponential mapping $\exp|U$ is $9$-determined in $h.$

$\exp|U$ is determined by the $9$-$jet$ of the metric
\end{theorem}

\begin{theorem}
\label{th6c3}

When there is collision between $Isoself_{c}$ and $Isoself_{1}$, then :

$-Isoself_{c}$ is the union of two curves.

$-Isoself_{1}$ is the union of one or two curves.

All the self-intersections are transversal and are not dependent on $\epsilon.$

There is a neighborhood $U$ of $\mathbb{S}\cap\{0<h<c\}$ (for small enough
values of $c$), such that :

The $7$-$jet$ of $CL$ : $CL^{7}=f_{4}h^{4}+f_{6}h^{6}+f_{7}h^{7}$ is a
sufficient jet of $CL$ on $U.$

The restriction of the exponential mapping $\exp|U$ is $7$-determined in $h.$

$\exp|U$ is determined by the $7$-$jet$ of the metric
\end{theorem}

\subsection{Proof of the stability results along the curve $C$}

\subsubsection{Definition and characterization of the self-intersections of
$CL$}

Let $Isoself$ denote the set of self-intersections of $CL$ on $C$, roughly
speaking the set of $(h,\varphi_{1},\varphi_{2})$ such that :

1) $\varphi_{1}\neq\varphi_{2},$

2) $CL(h,\varphi_{1})=CL(h,\varphi_{2}).$

We are interested in the germ of $Isoself$ along $\{h=0\}.$ A value of
$\varphi$ is said to be adherent to $Isoself$ if $(0,\varphi,\varphi
^{^{\prime}})$ lies in the closure of $Isoself$ for some $\varphi^{^{\prime}%
}.$ The set of some $\varphi$ is denoted by $A$-$Isoself.$

We know that if $\varphi\in A$-$Isoself$ then $\varphi^{^{\prime}}=\varphi
+\pi\in A$-$Isoself$ and there is not an other possibility (\cite{AEG}).

\begin{lemma}
\label{lem1} $A$-$Isoself\subset\{\varphi\,|\,\Psi(\varphi)=0\}.$
\end{lemma}

Assuming that $(h,\varphi,\,\varphi+\pi+\delta)\in A$-$Isoself,$ with $h>0,$
and $\delta$ small. Then :%

\begin{equation}
\overset{7}{\underset{i=4}{\sum}}f_{i}(\varphi+\pi+\delta)h^{i}=\overset
{7}{\underset{i=4}{\sum}}f_{i}(\varphi)h^{i}+O(h^{8}). \label{ef}%
\end{equation}

Since $f_{4}$ is independent of $\varphi$ and $f_{5}(\varphi)=0$ on $C$,
dividing \ref{ef} by $h^{6}$ we obtain :%

\[
f_{6}(\varphi+\delta)-f_{6}(\varphi)+h\,(f_{7}(\varphi+\delta+\pi
)-f_{7}(\varphi+\pi))+h\,(f_{7}(\varphi+\pi)-f_{7}(\varphi))+O(h^{2})=0
\]

Using $(B.P.)$ we obtain :
\[
\delta\frac{df_{6}}{d\varphi}(\varphi)-2hf_{7}(\varphi)+O(h^{2})+O(h\delta)=0
\]

Therefore if $\tilde{\varphi}\in A$-$Isoself,$ then $\frac{df_{6}}{d\varphi
}(\tilde{\varphi})\wedge f_{7}(\tilde{\varphi})=0\,$;\thinspace and lemma
\ref{lem1} is proved.

\begin{remark}
As a consequence of this lemma $A$-$Isoself$ splits into two subsets :
\[
A\text{-}Isoself=A\text{-}Isoself_{c}\cup A\text{-}Isoself_{1}.
\]

Where $A$-$Isoself_{c}=\{k\frac{\pi}{3}-\frac{\theta_{3}}{3}\}$ is the set of
cuspidal angles, and $A$-$Isoself_{1}$ is the set of the roots of $P(\varphi).$
\end{remark}

\subsubsection{Characterization of $Isoself_{c}$ on generic points of $C.$}

\smallskip\ 

We have to solve the equation :%

\[
\overset{7}{\underset{i=4}{\sum}}f_{i}(\varphi+\pi+\delta)h^{i}=\overset
{7}{\underset{i=4}{\sum}}f_{i}(\varphi)h^{i}+O(h^{8}),
\]

for small enough values of $\delta$ and $h.$

We know that this equation is equivalent to :

\smallskip\
\begin{equation}
f_{6}(\varphi+\delta)-f_{6}(\varphi)+h\,(f_{7}(\varphi+\delta+\pi
)-f_{7}(\varphi+\pi))+h\,(f_{7}(\varphi+\pi)-f_{7}(\varphi))+O(h^{2})=0,
\label{ec0}%
\end{equation}

For simplicities sake we denote$\dfrac{d}{d\varphi}$ by $^{\prime}$ ; hence
\ref{ec0} writes :%

\[
E_{c}=0=\delta\,f_{6}^{^{\prime}}(\varphi)+\frac{\delta^{2}}2f_{6}%
^{^{\prime\prime}}(\varphi)+\frac{\delta^{3}}6f_{6}^{^{\prime\prime\prime}%
}(\varphi)-h\,\delta\,f_{7}^{^{\prime}}(\varphi)-2h\,f_{7}(\varphi
)+O(\delta^{4})+O(h\,\delta^{2})+O(h^{2})
\]
By $(B.P.)$ $f_{6}^{^{\prime}}(0)=0$ (since $\theta_{3}=0$). Checking
$f_{6}^{^{\prime\prime\prime}}(\varphi)\wedge f_{6}^{^{\prime\prime}}%
(\varphi)$; we obtain :%

\[
f_{6}^{^{\prime\prime\prime}}(0)\wedge f_{6}^{^{\prime\prime}}%
(0)=-583200\left|  r_{3}\right|  \pi^{2}
\]

Hence :%

\[
E_{c}=0\Longleftrightarrow\left\{
\begin{array}
[c]{c}%
E_{c}\wedge f_{6}^{^{\prime\prime\prime}}(0)=0\\
\\
E_{c}\wedge f_{6}^{^{\prime\prime}}(0)=0
\end{array}
\right.
\]

Expanding $E_{c}$ we obtain%

\[%
\begin{tabular}
[c]{l}%
$E_{c}=\delta\,f_{6}^{^{\prime}}(0)+\delta\,\varphi\,f_{6}^{^{\prime\prime}%
}(0)+\frac\delta2\varphi^{2}f_{6}^{^{\prime\prime\prime}}(0)+O(\delta
\varphi^{3})+$\\
\\
$\frac{\delta^{2}}2\,f_{6}^{^{\prime\prime}}(0)+\frac{\delta^{2}}%
2\varphi\,f_{6}^{^{\prime\prime\prime}}(0)+O(\delta^{2}\varphi^{2}%
)+\frac{\delta^{3}}6f_{6}^{^{\prime\prime\prime}}(0)+O(\delta^{4}%
)-h\delta\,f_{7}^{^{\prime}}(0)+O(h\delta\varphi)$\\
\\
$-2hf_{7}(0)-2\,h\,\varphi\,f_{7}^{^{\prime}}(0)+O(h\varphi^{2})+O(h\,\delta
^{2})+O(h^{2})=0.$%
\end{tabular}
\]

Hence :%

\begin{equation}%
\begin{tabular}
[c]{l}%
$E_{c}\wedge f_{6}^{^{\prime\prime}}(0)=\frac\delta2\varphi^{2}f_{6}%
^{^{\prime\prime\prime}}(0)\wedge f_{6}^{^{\prime\prime}}(0)+\frac{\delta^{2}%
}2\varphi\,f_{6}^{^{\prime\prime\prime}}(0)\wedge f_{6}^{^{\prime\prime}%
}(0)+\frac{\delta^{3}}6f_{6}^{^{\prime\prime\prime}}(0)\wedge f_{6}%
^{^{\prime\prime}}(0)$\\
\\
$-2h\,f_{7}(0)\wedge f_{6}^{^{\prime\prime}}(0)+O(h\delta)+O(h\,\varphi
)+O(\delta\varphi^{3})+O(\delta^{4})+O(\varphi^{2}\,\delta^{2})+O(h^{2})=0.$%
\end{tabular}
\label{sec3}%
\end{equation}

and,
\begin{equation}%
\begin{array}
[c]{c}%
E_{c}\wedge f_{6}^{^{\prime\prime\prime}}(0)=(\delta\,\varphi\,+\frac
{\delta^{2}}2)f_{6}^{^{\prime\prime}}(0)\wedge f_{6}^{^{\prime\prime\prime}%
}(0)+\\
\\
O(h)+O(\delta\varphi^{2})+O(\delta^{3})+O(\delta^{2}\varphi)=0.
\end{array}
\label{sec4}%
\end{equation}

\smallskip\ We can solve \ref{sec3} in $h$ (implicit function theorem).
\begin{equation}
h=\frac{1}{2}\frac{f_{6}^{^{\prime\prime\prime}}(0)\wedge f_{6}^{^{\prime
\prime}}(0)}{\,f_{7}(0)\wedge f_{6}^{^{\prime\prime}}(0)}(\frac{\delta}%
{2}\varphi^{2}+\frac{\delta^{2}}{2}\varphi+\frac{\delta^{3}}{6})+\delta
^{2}\,O^{2}(\delta,\,\varphi)+O(\delta\,\varphi^{3}). \label{sec5}%
\end{equation}

We claim that $f_{7}(0)\wedge f_{6}^{^{\prime\prime}}(0)\neq0$.

In fact :
\[
\frac d{d\varphi}(f_{6}^{^{\prime}}\wedge f_{7})(0)=f_{6}^{^{\prime\prime}%
}(0)\wedge f_{7}(0),
\]

And zero is not a root of $P(\varphi)$ by theorem \ref{th1}.

\smallskip

Replacing $h$ by its value in \ref{sec4} we obtain :%

\begin{equation}
\left\{
\begin{array}
[c]{l}%
\varphi=-\frac\delta2+O(\delta^{2}),\\
\\
h=-\frac1{48}\frac{f_{6}^{^{\prime\prime\prime}}(0)\wedge f_{6}^{^{\prime
\prime}}(0)}{f_{7}(0)\wedge f_{6}^{^{\prime\prime}}(0)}\delta^{3}+O(\delta
^{4}).
\end{array}
\right.  \label{e6}%
\end{equation}

\begin{remark}
Since $h$ has to be positive, only one half of the curves defined by \ref{e6}
must be considered : either $\delta>0$ or $\delta<0$.
\end{remark}

\medskip\ \ 

Now let us show that this self-intersection is transversal. For this, the
following expression (the transversality-rate) $\mathcal{T}(\varphi,h,\delta)$
has to be nonzero on the self-intersection.%

\begin{equation}
\mathcal{T}(\varphi,h,\delta)=(f_{6}^{^{\prime}}(\varphi+\delta)+h\,f_{7}%
^{^{\prime}}(\varphi+\delta+\pi))\wedge(f_{6}^{^{\prime}}(\varphi
)+h\,f_{7}^{^{\prime}}(\varphi)+O(h^{2})). \label{e7}%
\end{equation}

Taking into account \ref{e6}, and expanding \ref{e7}, we obtain :
\begin{equation}
\mathcal{T}(\varphi,h,\delta)=-\frac{\delta^{3}}{24}f_{6}^{^{\prime
\prime\prime}}(0)\wedge f_{6}^{^{\prime\prime}}(0)+O(\delta^{4}) \label{e8}%
\end{equation}

Since,%

\[
-\frac{\delta^{3}}{24}f_{6}^{^{\prime\prime\prime}}(0)\wedge f_{6}%
^{^{\prime\prime}}(0)=2h\,\frac d{d\varphi}(f_{6}^{^{\prime}}\wedge
f_{7})(0)\neq0
\]

Hence \ref{e8} is nonzero.

\subsubsection{Characterization of $Isoself_{1}$ on generic points of $C$.}

Now we assume that $\theta_{3}\neq0$.

We want to solve again :
\begin{equation}
f_{6}(\varphi+\pi+\delta)+f_{7}(\varphi+\pi+\delta)\,h=f_{6}(\varphi
)+h\,f_{7}(\varphi)+O(h^{2}), \label{e10}%
\end{equation}
for $\varphi$ close to zero, zero being a simple root of $P(\varphi)$,
$\delta$ and $h$ small.

Expanding \ref{e10}, we obtain :%

\begin{equation}%
\begin{array}
[c]{l}%
0=\delta\,f_{6}^{^{\prime}}(\varphi)+O(\delta^{2})-h\,f_{7}(\varphi
+\delta)-h\,f_{7}(\varphi)+O(h^{2})\\
\\
0=\delta\,f_{6}^{^{\prime}}(\varphi)-2\,h\,f_{7}(\varphi)-h\,(\,f_{7}%
(\varphi+\delta)-\,f_{7}(\varphi))+O(h^{2})+O(\delta^{2}).
\end{array}
\label{e11}%
\end{equation}

This last expression rewrites :%

\begin{equation}%
\begin{array}
[c]{l}%
E_{1}=0=\delta\,f_{6}^{^{\prime}}(0)+\delta\,\varphi\,f_{6}^{^{\prime\prime}%
}(0)+O(\delta\varphi^{2})-2\,h\,f_{7}(0)-2\,h\,\varphi\,f_{7}^{^{\prime}%
}(0)+O(h^{2})+\\
\\
O(\delta^{2})+O(h\varphi^{2})+O(h\delta),
\end{array}
\label{e12}%
\end{equation}

Checking $\,f_{6}^{^{\prime}}(\varphi)\wedge f_{6}^{^{\prime\prime}}(\varphi)$
we obtain :%

\[
f_{6}^{^{\prime}}(\varphi)\wedge f_{6}^{^{\prime\prime}}\left(  \varphi
\right)  =32400\left|  r_{3}\right|  ^{2}\pi^{2}\,\sin^{2}(3\,\varphi
+\theta_{3}).
\]
Since $\theta_{3}$ has to be nonzero (far from cusp),%

\[
f_{6}^{^{\prime}}(0)\wedge f_{6}^{^{\prime\prime}}(0)\neq0.
\]

Therefore :%

\[
E_{1}=0\Longleftrightarrow\left\{
\begin{array}
[c]{c}%
E_{1}\wedge f_{6}^{^{\prime\prime}}(0)=0\\
\\
E_{1}\wedge f_{6}^{^{\prime}}(0)=0
\end{array}
\right.
\]

Hence :%

\begin{equation}%
\begin{array}
[c]{l}%
0=E_{1}\wedge f_{6}^{^{\prime\prime}}=\delta\,f_{6}^{^{\prime}}(0)\wedge
f_{6}^{^{\prime\prime}}(0)-2\,h\,f_{7}(0)\wedge f_{6}^{^{\prime\prime}%
}(0)+O(h^{2})\\
\\
+O(\delta^{2})+O(h\varphi)+O(h\delta)+O(\delta\varphi)
\end{array}
\label{se13}%
\end{equation}

and
\begin{equation}%
\begin{tabular}
[c]{l}%
$0=E_{1}\wedge f_{6}^{^{\prime}}=\delta\varphi\,f_{6}^{^{\prime\prime}%
}(0)\wedge f_{6}^{^{\prime}}(0)+\frac{\delta^{2}}2\,f_{6}^{^{\prime\prime}%
}(0)\wedge f_{6}^{^{\prime}}(0)$\\
\\
$-2\,h\,f_{7}(0)\wedge f_{6}^{^{\prime}}(0)-2\,h\varphi\,f_{7}^{^{\prime}%
}(0)\wedge f_{6}^{^{\prime}}(0)-\,h\delta\,f_{7}^{^{\prime}}(0)\wedge
f_{6}^{^{\prime}}(0)+$\\
\\
$O(\delta^{2}\varphi)+O(\delta^{3})+O(h\varphi^{2})+O(h\,\delta^{2}%
)+O(\delta\,\varphi^{2})+O(h^{2}\,\delta)+O(h^{3})$%
\end{tabular}
\label{se14}%
\end{equation}

\smallskip\ 

If $\,f_{7}(0)\wedge f_{6}^{^{\prime\prime}}(0)\neq0,(i.e.$ $f_{7}(0)\neq0$
and $f_{7}(0)=\lambda f_{6}^{^{\prime}}(0)),$ then $\ref{se13}$ gives us :
\[
\delta=2\frac{f_{7}(0)\wedge f_{6}^{^{\prime\prime}}(0)}{f_{6}^{^{\prime}%
}(0)\wedge f_{6}^{^{\prime\prime}}(0)}\,h+h\,O^{1}(\varphi,\,h)=2\lambda
\,h+h\,O^{1}(\varphi,\,h),\;\lambda\neq0
\]

Replacing $\delta$ by its value in \ref{se14} and dividing by $h$ ($h>0$) ;
$E_{1}\wedge f_{6}^{^{\prime}}(0)=0$ become :
\begin{equation}
2\varphi\frac{d}{d\varphi}(f_{6}^{^{\prime}}\wedge f_{7})(0)-2\,f_{7}(0)\wedge
f_{6}^{^{\prime}}(0)+2\lambda\,h\frac{d}{d\varphi}(f_{6}^{^{\prime}}\wedge
f_{7})(0)+O(h^{2})+O(\varphi^{2})=0 \label{sepp}%
\end{equation}
Recalling that $f_{7}(0)\wedge f_{6}^{^{\prime}}(0)=0$, and solving \ref{sepp}
in $\varphi$; we obtain at the end :
\begin{equation}
\left\{
\begin{array}
[c]{l}%
\varphi=-\lambda h+h\,O(h)\\
\\
\delta=2\lambda h+h\,O(h)
\end{array}
\right.  \label{se1f}%
\end{equation}

Since $h>0$ the curve defined by \ref{se1f} is a smooth one starting from zero.

\medskip\ 

We have to show again that this self-intersection is transversal.

The transversality-rate $\mathcal{T}_{1}(\varphi,h,\delta)$ has to be nonzero
on the self-intersection.
\[
\mathcal{T}_{1}(\varphi,h,\delta)=(f_{6}^{^{\prime}}(\varphi+\delta
)+h\,f_{7}^{^{\prime}}(\varphi+\delta+\pi))\wedge(f_{6}^{^{\prime}}%
(\varphi)+h\,f_{7}^{^{\prime}}(\varphi)+O(h^{2})).
\]%

\[
\mathcal{T}_{1}(\varphi,h,\delta)=(f_{6}^{^{\prime}}(\varphi)-h\,f_{7}%
^{^{\prime}}(\varphi)+\delta\,f_{6}^{^{\prime\prime}}(\varphi)+O(\delta
^{2})+O(\delta\,h))\wedge(f_{6}^{^{\prime}}(\varphi)+h\,f_{7}^{^{\prime}%
}(\varphi)+O(h^{2})).
\]

It is easy to see that :
\begin{equation}
\mathcal{T}_{1}(\varphi,h,\delta)=h\,\frac d{d\varphi}(f_{6}^{^{\prime}}\wedge
f_{7})+O(h^{2})=O(h^{2}). \label{t1}%
\end{equation}

Since zero is not a double root of $f_{6}^{^{\prime}}\wedge f_{7},$ then
$\mathcal{T}(\varphi,h,\delta)\neq0$ for $\varphi$ close to zero.

From now we assume that $f_{7}(0)=0.$

Hence :
\[
\frac d{d\varphi}(f_{6}^{^{\prime}}\wedge f_{7})(0)=f_{6}^{^{\prime}}(0)\wedge
f_{7}^{^{\prime}}(0)\neq0.
\]

In this case we will need higher order jet of $CL$ : namely $f_{9}(\varphi)$.

\begin{remark}
$f_{8}$ is insufficient because of the fact that $f_{8}(\varphi+\pi
)=f_{8}(\varphi).$
\end{remark}

We have to solve now :%

\begin{equation}%
\begin{tabular}
[c]{l}%
$E_{1}^{^{\prime}}=\delta\,f_{6}^{^{\prime}}(\varphi)+\frac{\delta^{2}%
}2\,f_{6}^{^{\prime\prime}}(\varphi)+\frac{\delta^{3}}6\,f_{6}^{^{\prime
\prime}}(\varphi)+O(\delta^{4})-2h\,f_{7}(\varphi)-h\delta\,f_{7}^{^{\prime}%
}(\varphi)+$\\
\\
$h^{2}\delta\,f_{8}^{^{\prime}}(0)+O(h^{3}\delta)+O(h\delta^{2})-2h^{3}%
f_{9}(0)+O(h^{4})=0.$%
\end{tabular}
\label{e1p0}%
\end{equation}%

\[
f_{6}^{^{\prime}}(0)\wedge f_{7}^{^{\prime}}(0)\neq0.
\]

Hence :%

\[
E_{1}^{^{\prime}}=0\Longleftrightarrow\left\{
\begin{array}
[c]{c}%
E_{1}^{^{\prime}}\wedge f_{7}^{^{\prime}}(0)=0\\
\\
E_{1}^{^{\prime}}\wedge f_{6}^{^{\prime}}(0)=0
\end{array}
\right.
\]

Also :
\begin{equation}%
\begin{tabular}
[c]{l}%
$0=E_{1}^{^{\prime}}\wedge f_{7}^{^{\prime}}=\delta\,f_{6}^{^{\prime}%
}(0)\wedge f_{7}^{^{\prime}}(0)-2\,h^{3}\,f_{9}(0)\wedge f_{7}^{^{\prime}%
}(0)+O(h^{4})+O(\delta^{2})+$\\
\\
$O(\varphi\,h)+O(h\,\delta)+O(\delta\,\varphi).$%
\end{tabular}
\label{se1p3}%
\end{equation}

and%

\begin{equation}%
\begin{tabular}
[c]{l}%
$0=E_{1}^{^{\prime}}\wedge f_{6}^{^{\prime}}=-2h\varphi\,f_{7}^{^{\prime}%
}\wedge f_{6}^{^{\prime}}-2\,h\,^{3}f_{9}(0)\wedge f_{6}^{^{\prime}}(0)$\\
\\
$O(h\varphi^{2})+O(\delta^{2})+O(\delta\,\varphi)+O(h\,\delta)+O(h^{4})$%
\end{tabular}
\label{se1p4}%
\end{equation}

\smallskip\ 

From \ref{se1p3} one obtains :
\[
\delta=2\frac{f_{9}(0)\wedge f_{7}^{^{\prime}}(0)}{f_{6}^{^{\prime}}(0)\wedge
f_{7}^{^{\prime}}(0)}\,h^{3}+h^{3}\,O^{1}(\varphi,\,h)
\]

Replacing $\delta$ by its value in \ref{se1p4}, and dividing by $h$ ($h>0$),
$E_{1}^{^{\prime}}\wedge f_{6}^{^{\prime}}(0)=0$ become :
\[
-2\varphi\,f_{7}^{^{\prime}}(0)\wedge f_{6}^{^{\prime}}(0)-2\,h^{2}%
f_{9}(0)\wedge f_{6}^{^{\prime}}(0)+O(h^{3})+O(\varphi^{2})=0
\]

Hence :
\begin{equation}
\left\{
\begin{array}
[c]{l}%
\varphi=A^{^{\prime}}h^{2}+h\,O^{2}(h)\\
\\
\delta=B^{^{\prime}}\,h^{3}+h\,O^{3}(h)
\end{array}
\right.  \label{se1pf}%
\end{equation}

Again the curve defined by \ref{se1pf} is a smooth one starting from zero.

\medskip\ 

Showing now that this self-intersection is transversal.

\smallskip\ 

The transversality-rate $\mathcal{T}_{2}(\varphi,h,\delta)$ has to be nonzero
on the self-intersection.
\[
\mathcal{T}_{2}(\varphi,h,\delta)=(f_{6}^{^{\prime}}(\varphi+\delta
)+h\,f_{7}^{^{\prime}}(\varphi+\delta+\pi))\wedge(f_{6}^{^{\prime}}%
(\varphi)+h\,f_{7}^{^{\prime}}(\varphi)+O(h^{2})).
\]%

\[
\mathcal{T}_{2}(\varphi,h,\delta)=(f_{6}^{^{\prime}}(\varphi)-h\,f_{7}%
^{^{\prime}}(\varphi)+O(\delta)+O(\delta\,h))\wedge(f_{6}^{^{\prime}}%
(\varphi)+h\,f_{7}^{^{\prime}}(\varphi)+O(h^{2})).
\]
\smallskip

As in the previous case :
\begin{equation}
\mathcal{T}_{2}(\varphi,h,\delta)=h\,\frac d{d\varphi}(f_{6}^{^{\prime}}\wedge
f_{7})+O(h^{2})=O(h^{2}). \label{t2}%
\end{equation}

Zero is not a double root of $f_{6}^{^{\prime}}\wedge f_{7},$ hence
$\mathcal{T}_{2}(\varphi,h,\delta)\neq0$ close to zero.

\subsubsection{Characterization de $Isoself_{1}$ on isolated points of $C$}

\paragraph{Double roots of P($\varphi$).\newline }

We recall hypothesis :

\smallskip\ 

i) $\frac d{d\varphi}(f_{6}^{^{\prime}}\wedge f_{7})(0)=f_{6}^{^{\prime\prime
}}(0)\wedge f_{7}(0)+f_{6}^{\prime}(0)\wedge f_{7}^{^{\prime}}(0)=0.$

ii) $\theta_{3}\neq0.$

The equation of the self-intersection is the following :%

\begin{equation}%
\begin{tabular}
[c]{l}%
$E_{1}^{1}=\delta\,f_{6}^{^{\prime}}(\varphi)+\frac{\delta^{2}}2\,f_{6}%
^{^{\prime\prime}}(\varphi)+\frac{\delta^{3}}6\,f_{6}^{^{\prime\prime}%
}(\varphi)+O(\delta^{4})-2h\,f_{7}(\varphi)-h\delta\,f_{7}^{^{\prime}}%
(\varphi)+$\\
\\
$h^{2}\delta\,f_{8}^{^{\prime}}(0)+O(h^{3}\delta)+O(h\delta^{2})-2h^{3}%
f_{9}(0)+O(h^{4})=0$%
\end{tabular}
\label{se110}%
\end{equation}

Our hypothesis ensure us that :
\[
f_{6}^{^{\prime}}(0)\wedge f_{6}^{^{\prime\prime}}(0)\neq0.
\]

Hence :%

\[
E_{1}^{1}=0\Longleftrightarrow\left\{
\begin{array}
[c]{c}%
E_{1}^{1}\wedge f_{6}^{^{\prime\prime}}(0)=0\\
\\
E_{1}^{1}\wedge f_{6}^{^{\prime}}(0)=0
\end{array}
\right.
\]

We interesse at first to $E_{1}^{1}\wedge f_{6}^{^{\prime\prime}}(0)=0.$

We obtain :
\[%
\begin{tabular}
[c]{l}%
$0=\delta f_{6}^{^{\prime}}(0)\wedge f_{6}^{^{\prime\prime}}(0)-2hf_{7}%
(0)\wedge f_{6}^{^{\prime\prime}}(0)+$\\
\\
$O(\delta^{2})+O(\delta\varphi)+O(h^{2})+O(h\varphi)+O(h\delta)$%
\end{tabular}
\]

Solving this last equation in $\delta$ ; one obtains :%

\begin{equation}
\delta=2\frac{f_{7}(0)\wedge f_{6}^{^{\prime\prime}}(0)}{f_{6}^{^{\prime}%
}(0)\wedge f_{6}^{^{\prime\prime}}(0)}\,h+hO^{1}(h,\varphi). \label{g12}%
\end{equation}

We will distinguish two cases.

\smallskip\ 

1) $\,f_{7}(0)\neq0.$

In that case :%

\begin{equation}%
\begin{tabular}
[c]{l}%
$f_{7}(0)=\lambda\,f_{6}^{^{\prime}}(0)$ with $\lambda\neq0.$\\
and\\
$\delta=2\lambda\,h+O(h^{2})+O(h\varphi).$%
\end{tabular}
\label{l1}%
\end{equation}

We will consider now $E_{1}^{1}\wedge f_{6}^{^{\prime}}(0)=0$. Assume that
there is at least one solution.

This solution is on the form:
\[
\varphi=\chi\,h+\chi_{2}\,h^{2}+\chi_{2}\,h^{3}+O(h^{4}).
\]

It is easy to see that $\chi\,\neq0.$

Replacing $\varphi$ by its value and dividing by $h^{3}$ we obtain :%

\[%
\begin{tabular}
[c]{l}%
$m\,\chi^{2}+2\lambda m\chi+n+O(h)=0.$\\
\\
where,\\
$m=\lambda\,f_{6}^{^{\prime\prime\prime}}(0)\wedge f_{6}^{^{\prime}}%
(0)+f_{6}^{^{\prime}}(0)\wedge f_{7}^{^{^{\prime\prime}}}(0)$\\
\\
$n=\frac{4\lambda}3f_{6}^{^{\prime\prime\prime}}(0)\wedge f_{6}^{^{\prime}%
}(0)-2\lambda f_{8}^{^{\prime}}(0)\wedge f_{6}^{^{\prime}}(0)-2\lambda
f_{9}(0)\wedge f_{6}^{^{\prime}}(0)$%
\end{tabular}
\]

By remark \ref{RC} we have the following important fact : $n\neq0.$

Hence : $m\neq0.$

Denoting $\chi_{1}$ and $\chi_{2}$ two presumed solutions of
\begin{equation}
m\,\chi^{2}+2\lambda m\chi+n=0. \label{eq2}%
\end{equation}
Then :%

\begin{equation}
\chi_{1}+\chi_{2}=-2\lambda\label{som1}%
\end{equation}

\smallskip\ 

If $\chi_{1}=\chi_{2}=\chi_{0}=-\lambda,$ then $-m\lambda^{2}+n=0$. This is
impossible by remark \ref{RC} once more.

Therefore \ref{eq2} has not a double root on $\mathbb{R}$.

Hence in this case there is either no self-intersection or two
self-intersections which satisfy the following equations :%

\[%
\begin{tabular}
[c]{l}%
$\delta=2\lambda h+O(h^{2})$\\
\\
$\varphi=\chi h+O(h^{2})$ where, $\chi=\chi_{i}$ et $\,i=1,2$%
\end{tabular}
\]

It remains to show that theses self-intersections are transversal.

Denoting by $\mathcal{T}_{1}^{1}(\varphi,h,\delta)$ the transversality-rate,
we have as usually
\[
\mathcal{T}_{1}^{1}(\varphi,h,\delta)=(f_{6}^{^{\prime}}(\varphi
+\delta)+h\,f_{7}^{^{\prime}}(\varphi+\delta+\pi))\wedge(f_{6}^{^{\prime}%
}(\varphi)+h\,f_{7}^{^{\prime}}(\varphi)+O(h^{2})).
\]

Expanding $\mathcal{T}_{1}^{1}(\varphi,h,\delta)$, we obtain :
\begin{equation}
\mathcal{T}_{1}^{1}(\varphi,h,\delta)=2h^{2}(\chi+\lambda)\frac{d^{2}%
}{d\varphi^{2}}(f_{6}^{^{\prime}}\wedge f_{7})(0)+O(h^{3}). \label{t11}%
\end{equation}

Since $(\chi+\lambda)\neq0,$ then $\mathcal{T}_{1}^{1}(\varphi,h,\delta)\neq0$
on our two self-intersections, hence they are well transversal.

\smallskip\ 

2) $f_{7}(0)=0$

Now $\lambda=0$ and $f_{6}^{^{\prime}}(0)\wedge f_{7}^{^{\prime}}(0)=0.$

Set $f_{7}^{^{\prime}}(0)=\gamma\,f_{6}^{^{\prime}}(0)$

\begin{remark}
Since zero is not a double root of $P(\varphi),$then $\gamma\neq0.$
\end{remark}

Taking again $E_{1}^{1}\wedge f_{6}^{^{^{\prime\prime}}}(0),i.e.$ :%

\[%
\begin{tabular}
[c]{l}%
$0=\delta f_{6}^{^{\prime}}(0)\wedge f_{6}^{^{\prime\prime}}(0)-2h\varphi
f_{7}^{^{\prime}}(0)\wedge f_{6}^{^{\prime\prime}}(0)+$\\
\\
$O(\delta^{2})+O(\delta\varphi^{2})+O(h^{3})+O(h\varphi^{2})+O(h\delta).$%
\end{tabular}
\]

We obtain :%

\[
\delta=2\gamma h\varphi+(h\varphi^{2})+O(h^{3})
\]

Replacing once more $\delta$ by its value in $E_{1}^{1}\wedge f_{6}^{^{\prime
}}(0)=0,$with the hypothesis :
\[
\varphi=\chi\,h+O(h^{2}).
\]

Dividing by $h^{3}$ we obtain :%

\[%
\begin{tabular}
[c]{l}%
$m\chi^{2}-2f_{9}(0)\wedge f_{6}^{^{\prime}}(0)+O(h)=0$\\
\\
where\\
\\
$m=2f_{6}^{^{\prime\prime}}(0)\wedge f_{7}^{^{\prime}}(0)+f_{6}^{^{\prime}%
}(0)\wedge f_{7}^{^{\prime\prime}}(0)=\frac{d^{2}}{d\varphi^{2}}%
(f_{6}^{^{\prime}}\wedge f_{7})(0)\neq0.$%
\end{tabular}
\]

Thus ; $\chi=\pm\sqrt{\kappa}$ where $\kappa=\frac{2f_{9}(0)\wedge
f_{6}^{^{\prime}}(0)}m\neq0.$

\smallskip\ 

Therefore if $\kappa<0$ there is no self-intersection and if $\kappa>0,$there
two self-intersections which satisfy the following equations :%

\[%
\begin{tabular}
[c]{l}%
$\delta=2\gamma\chi h^{2}+O(h^{2})$\\
\\
$\varphi=\pm\chi h+O(h^{2})$%
\end{tabular}
\]

As usually we will show that these two self-intersections are transversal .

Denoting by $\mathcal{T}_{1}^{2}(\varphi,h,\delta)$ the tranversality-rate ;
we have :
\[
\mathcal{T}_{1}^{2}(\varphi,h,\delta)=(f_{6}^{^{\prime}}(\varphi
+\delta)+h\,f_{7}^{^{\prime}}(\varphi+\delta+\pi))\wedge(f_{6}^{^{\prime}%
}(\varphi)+h\,f_{7}^{^{\prime}}(\varphi)+O(h^{2})).
\]

Expanding this last expression, we obtain :
\begin{equation}
\mathcal{T}_{1}^{2}(\varphi,h,\delta)=\pm2\chi h^{2}\frac{d^{2}}{d\varphi^{2}%
}(f_{6}^{^{\prime}}\wedge f_{7})(0)+O(h^{3}). \label{t11b}%
\end{equation}

Taking account our hypothesis, we know that :
\[
2\chi h^{2}\frac{d^{2}}{d\varphi^{2}}(f_{6}^{^{\prime}}\wedge f_{7})(0)\neq0.
\]

\subparagraph{Collision between $isoself_{1}$ and $isoself_{c}$.\newline }

Now, zero is a cusp point and a simple root of $P(\varphi)$.

Hence
\[
f_{6}^{^{\prime}}(0)=0\text{ et }\frac d{d\varphi}(f_{6}^{^{\prime}}\wedge
f_{7})(0)=f_{6}^{^{\prime\prime}}(0)\wedge f_{7}(0)=0.
\]

In compensation
\begin{equation}
\frac{d^{2}}{d\varphi^{2}}(f_{6}^{^{\prime}}\wedge f_{7})(0)=f_{6}%
^{^{\prime\prime\prime}}(0)\wedge f_{7}(0)+2\,f_{6}^{^{\prime\prime}}(0)\wedge
f_{7}^{^{\prime}}(0)\neq0. \label{cl}%
\end{equation}

and%

\[
f_{6}^{^{\prime\prime\prime}}(0)\wedge f_{6}^{^{\prime\prime}}(0)\neq0.
\]

Again, we will distinguish two cases.

\smallskip\ 

1)$f_{7}(0)\neq0.$

In this case, there exist $\lambda\neq0$ such that $f_{7}(0)=\lambda
f_{6}^{^{\prime\prime}}(0)$

The equation of the self-intersection is :%

\begin{equation}%
\begin{tabular}
[c]{r}%
$E_{c}^{1}=\delta\,f_{6}^{^{\prime}}(\varphi)+\frac{\delta^{2}}2f_{6}%
^{^{\prime\prime}}(\varphi)+\frac{\delta^{3}}6f_{6}^{^{\prime\prime\prime}%
}(\varphi)-h\,\delta\,f_{7}^{^{\prime}}(\varphi)-2h\,f_{7}(\varphi)+$\\
$O(\delta^{4})+O(h\,\delta^{2})+O(h^{2})=0.$%
\end{tabular}
\label{esc1}%
\end{equation}%

\[
E_{c}^{1}=0\Longleftrightarrow\left\{
\begin{array}
[c]{c}%
E_{c}^{1}\wedge f_{6}^{^{\prime\prime\prime}}(0)=0\\
\\
E_{c}^{1}\wedge f_{6}^{^{\prime\prime}}(0)=0
\end{array}
\right.
\]%

\[%
\begin{tabular}
[c]{l}%
$E_{c}^{1}=\delta\,f_{6}^{^{\prime}}(0)+\delta\,\varphi\,f_{6}^{^{\prime
\prime}}(0)+\frac\delta2\varphi^{2}f_{6}^{^{\prime\prime\prime}}%
(0)+O(\delta\varphi^{3})+$\\
\\
$\frac{\delta^{2}}2\,f_{6}^{^{\prime\prime}}(0)+\frac{\delta^{2}}%
2\varphi\,f_{6}^{^{\prime\prime\prime}}(0)+O(\delta^{2}\varphi^{2}%
)+\frac{\delta^{3}}6f_{6}^{^{\prime\prime\prime}}(0)+O(\delta^{4}%
)-h\delta\,f_{7}^{^{\prime}}(0)+O(h\delta\varphi)$\\
\\
$-2hf_{7}(0)-2\,h\,\varphi\,f_{7}^{^{\prime}}(0)+O(h\varphi^{2})+O(h\,\delta
^{2})+O(h^{2})=0.$%
\end{tabular}
\]

$f_{6}^{^{\prime}}(0)=0,$ Therefore :%

\begin{equation}%
\begin{tabular}
[c]{l}%
$E_{c}^{1}\wedge f_{6}^{^{\prime\prime}}(0)=\frac\delta2\varphi^{2}%
f_{6}^{^{\prime\prime\prime}}(0)\wedge f_{6}^{^{\prime\prime}}(0)+\frac
{\delta^{2}}2\varphi\,f_{6}^{^{\prime\prime\prime}}(0)\wedge f_{6}%
^{^{\prime\prime}}(0)+\frac{\delta^{3}}6f_{6}^{^{\prime\prime\prime}}(0)\wedge
f_{6}^{^{\prime\prime}}(0)$\\
\\
$-2h\,f_{7}(0)\wedge f_{6}^{^{\prime\prime}}(0)-2h(\varphi+\frac\delta
2)f_{7}^{^{\prime}}(0)\wedge f_{6}^{^{\prime\prime}}(0)+$\\
\\
$O(h\,\varphi^{2})+O(\delta\varphi^{3})+O(\delta^{4})+O(\varphi^{2}%
\,\delta^{2})+O(h^{2})=0,$%
\end{tabular}
\label{ec12}%
\end{equation}

and%

\begin{equation}%
\begin{array}
[c]{c}%
E_{c}^{1}\wedge f_{6}^{^{\prime\prime\prime}}(0)=(\delta\,\varphi
\,+\frac{\delta^{2}}2-2\lambda h)f_{6}^{^{\prime\prime}}(0)\wedge
f_{6}^{^{\prime\prime\prime}}(0)\\
\\
O(h\varphi)+O(\delta\varphi^{2})+O(\delta^{3})+O(h^{2})+O(h\delta)=0.
\end{array}
\label{ec13}%
\end{equation}

We can solve \ref{ec13} in $h.$%

\begin{equation}
h=\frac\delta{2\lambda}(\varphi+\frac\delta2)+O(\delta\varphi^{2}%
)+O(\delta^{3}) \label{vh}%
\end{equation}

Taking account the previous characterization of $isoself_{c}$ and
$isoself_{1},$ one can assume that :
\[
\varphi=\frac{-\delta}2+O(\delta^{2}).
\]

\begin{remark}
If $\varphi\neq\frac{-\delta}{2}$, then $h=h(\delta^{2})+O(\delta^{3}).$ It is
easy to see that this case is not generic.
\end{remark}

Replacing $\varphi$ by $\frac{-\delta}{2}$ in \ref{vh} we obtain
$h=h(\delta^{3})$.

The equation \ref{ec12} become :%

\begin{equation}
\frac{\delta^{2}}{24}f_{6}^{^{\prime\prime}}(0)\wedge f_{6}^{^{\prime
\prime\prime}}(0)+O^{3}(\delta)=0. \label{fic}%
\end{equation}

Which is not possible because $f_{6}^{^{\prime\prime}}(0)\wedge f_{6}%
^{^{\prime\prime\prime}}(0)\neq0.$

Hence in that case there is no self-intersection.

\subsection{The stability}

In this case the suspension of the exponential mapping is not a ``Whitney
map''. The main argument here is the Mather theorem.

If $\beta_{2,2}=0$, then $\beta_{3,3}$ $($or $r_{3})$ and $\mu$ are nonzero.
The intersection $CL_{w}^{\pm}$ between $CL$ and the planes $w=c$ (for small
enough values of $\left|  w\right|  $ ) is a closed curve having fold-points,
cusp-points (six) and self-intersections (see figures).

According to our previous subsection ``characterization of self-intersection
on $C$'' (Which finds here its justification), all the generic
self-intersections of $CL_{w}^{\pm}$ are transversal.

Thus $CL$ is in ``general position '', which ends the proofs of theorems
\ref{th6c1}, \ref{th6c2} and \ref{th6c3}.

\paragraph{Figures.}%

\begin{center}
\includegraphics[
height=3.9755in,
width=5.3108in
]%
{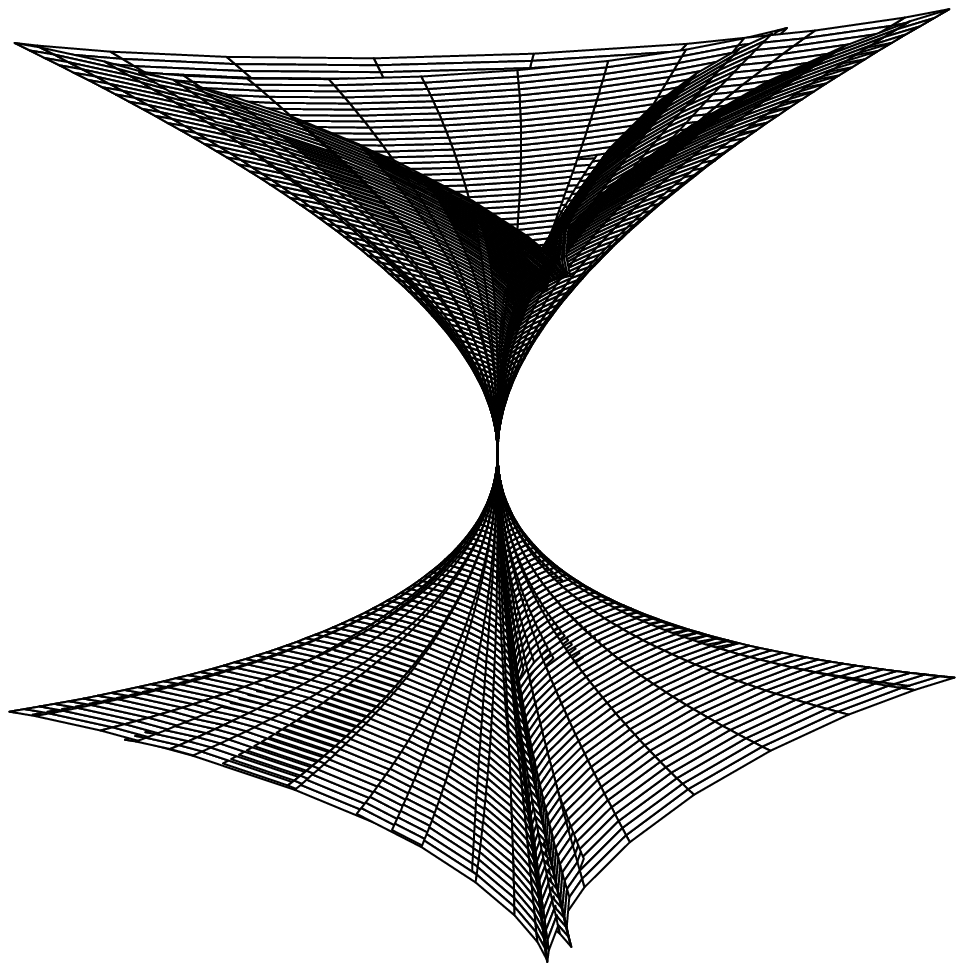}%
\\
$CL^{7}:$ when $r_{2}=0$ and $r_{3}\neq0$%
\end{center}
\begin{center}
\includegraphics[
height=3.9755in,
width=5.3108in
]%
{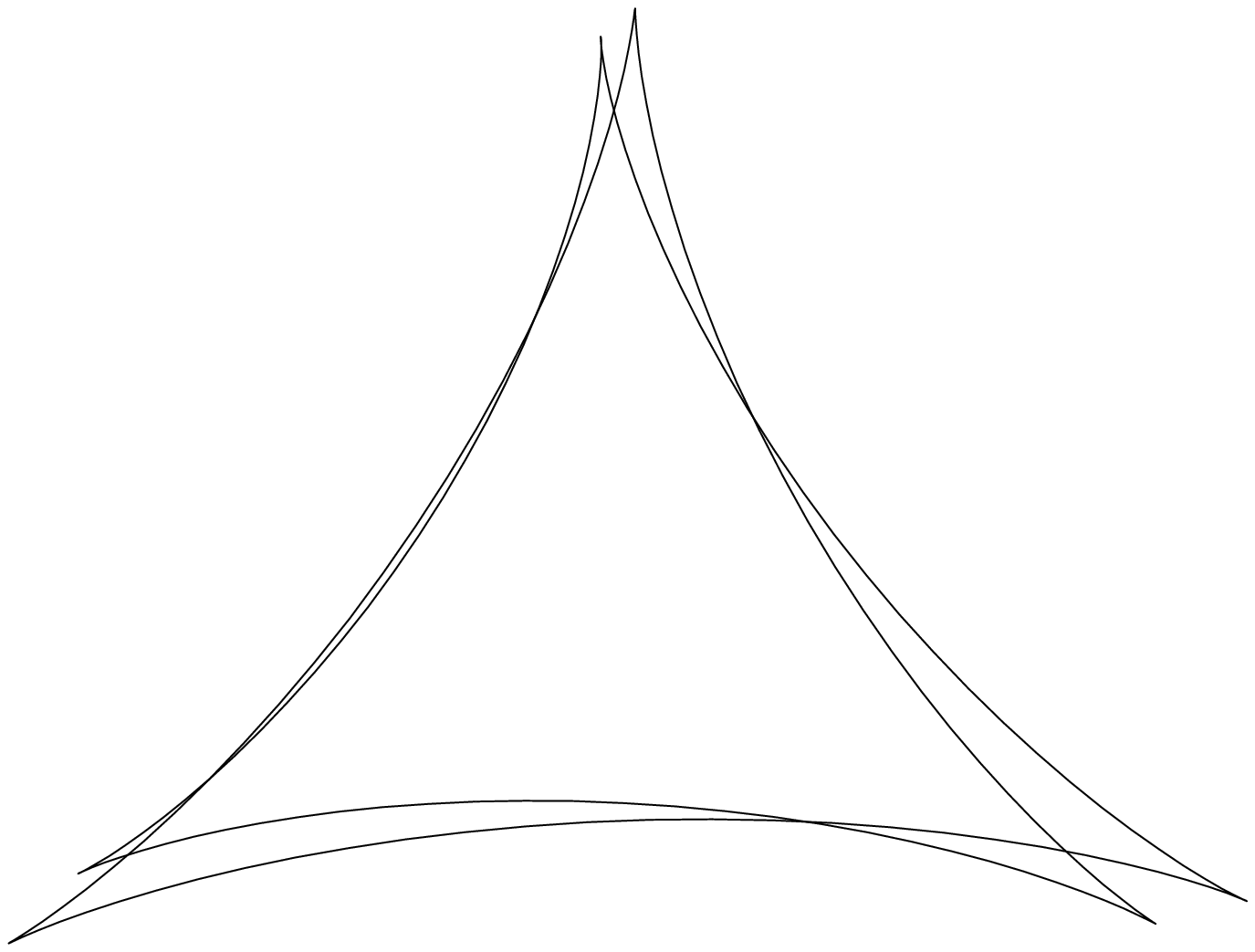}%
\\
$CL_{w}^{\pm}\,$\ (an example with $7$-$isoself$)
\end{center}
\begin{center}
\includegraphics[
height=3.9755in,
width=5.3108in
]%
{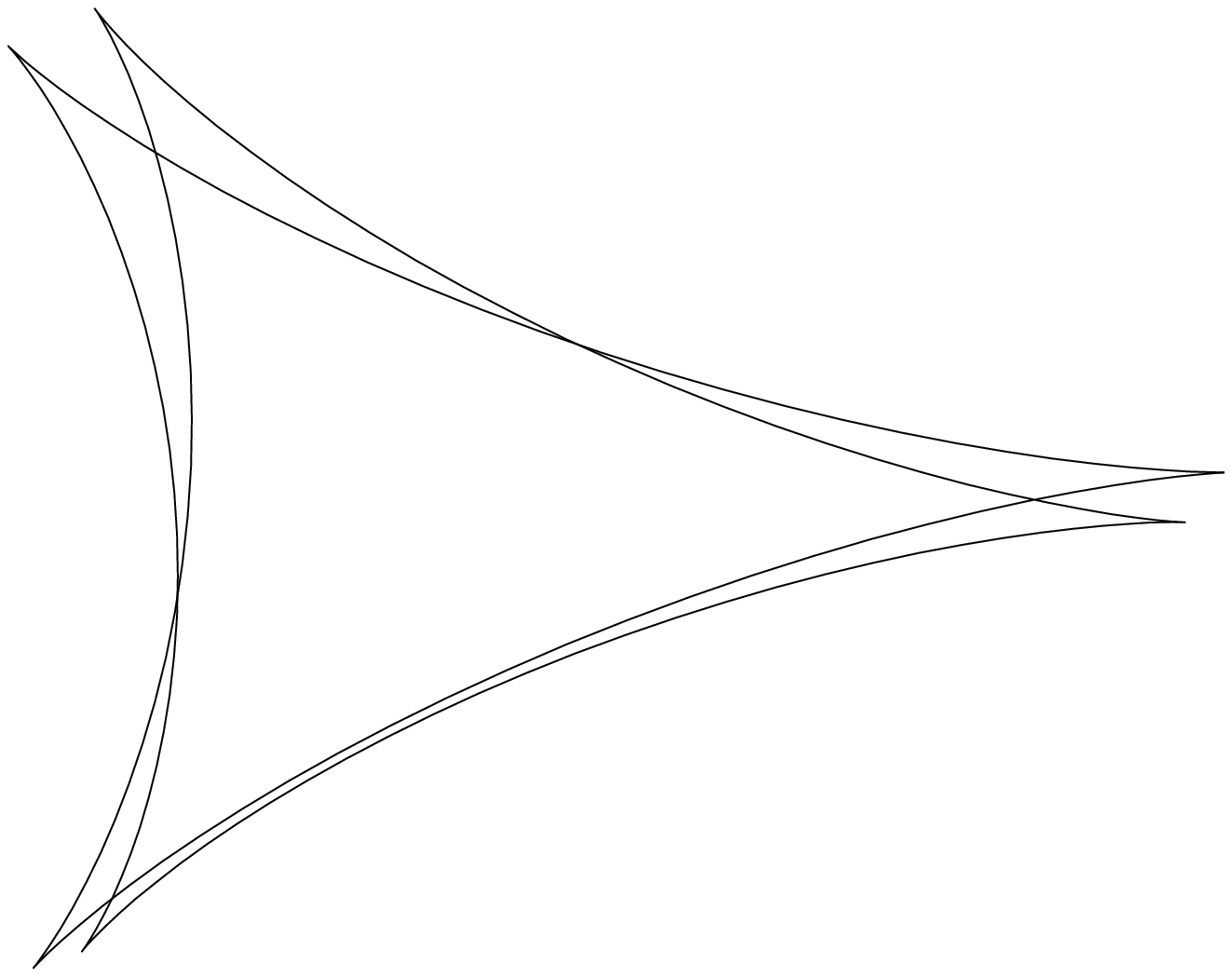}%
\\
$CL_{w}^{\pm}\,$\ (an example with $5$-$isoself$)
\label{CLP2}%
\end{center}

\newpage

\section{Appendix}

\subsection{Figures: wave front and conjugate loci\label{appA0}}%

\begin{center}
\includegraphics[
trim=0.000000in 0.000000in 0.002402in 0.002393in,
height=10.0957cm,
width=13.4873cm
]%
{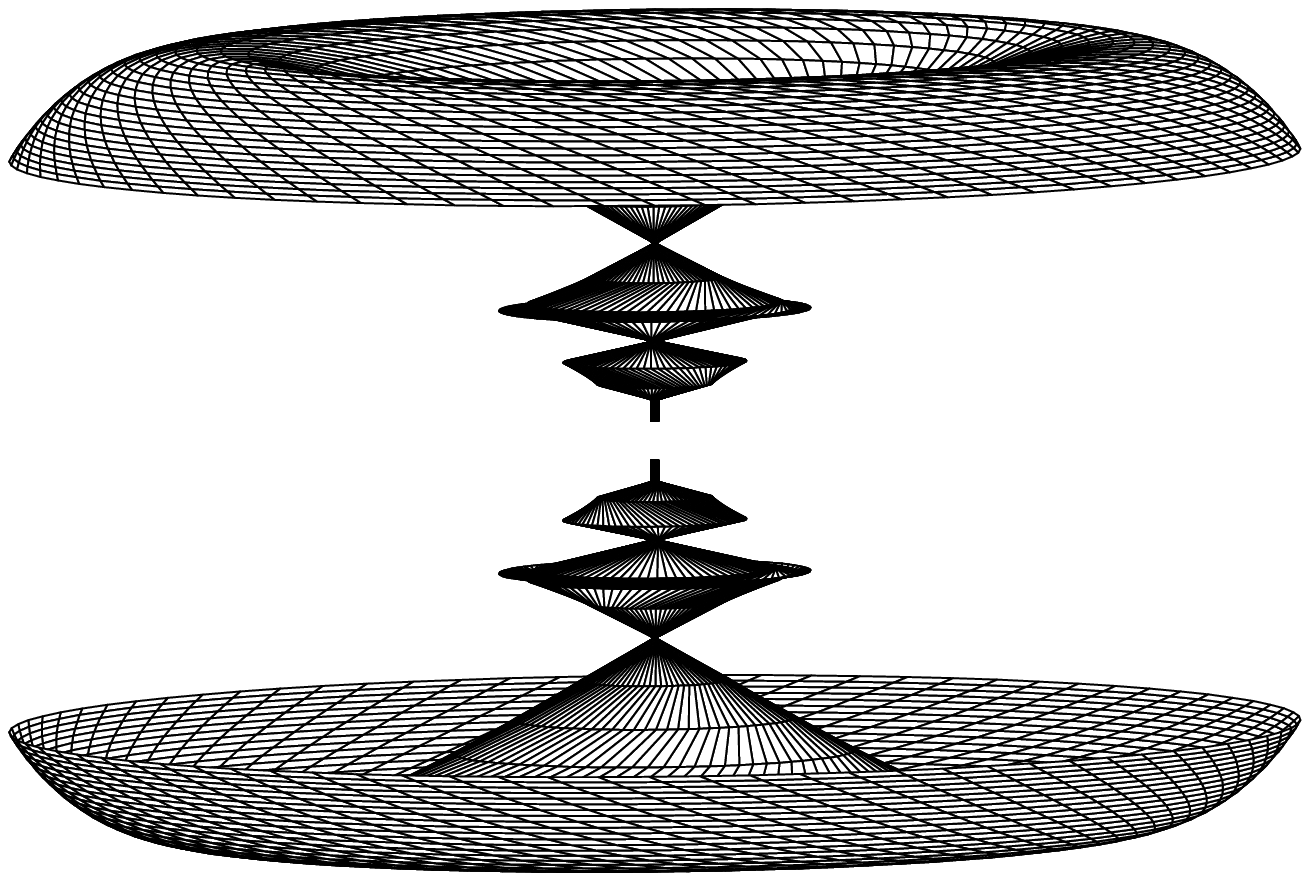}%
\\
The waves front when $r_{2}\neq0$ .
\label{chap2f1}%
\end{center}

\begin{center}
\includegraphics[
trim=0.000000in 0.000000in 0.002402in 0.002393in,
height=10.0957cm,
width=13.4873cm
]%
{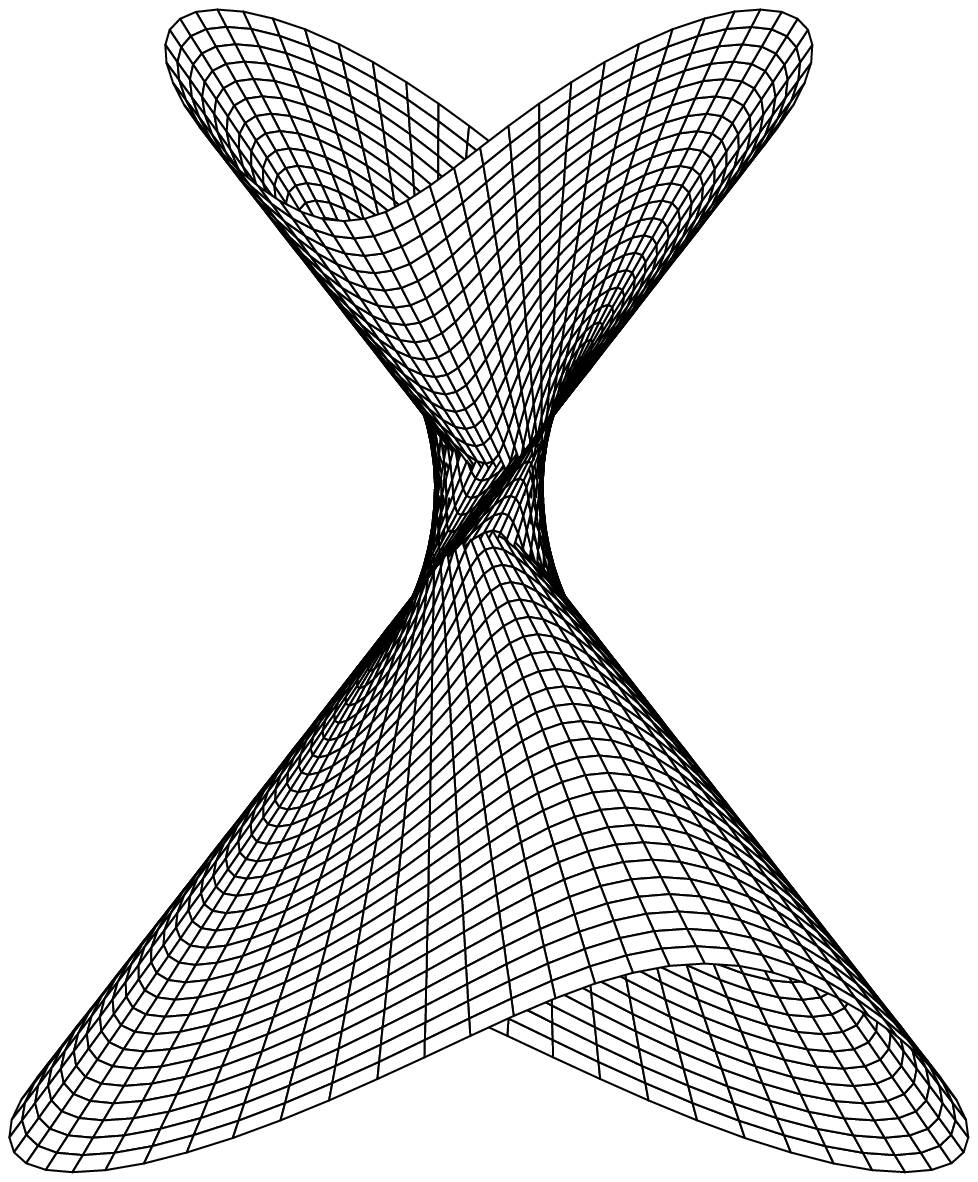}%
\\
Zoom on the zone C
\label{chap2f2}%
\end{center}

\newpage%

\begin{center}
\includegraphics[
trim=0.000000in 0.000000in 0.002402in 0.002393in,
height=10.0957cm,
width=13.4873cm
]%
{ffob1.1.eps}%
\\
The waves front when $r_{2}=0$ and $r_{3}\neq0$.
\label{chapt2f3}%
\end{center}

\begin{center}
\includegraphics[
height=3.9755in,
width=5.3117in
]%
{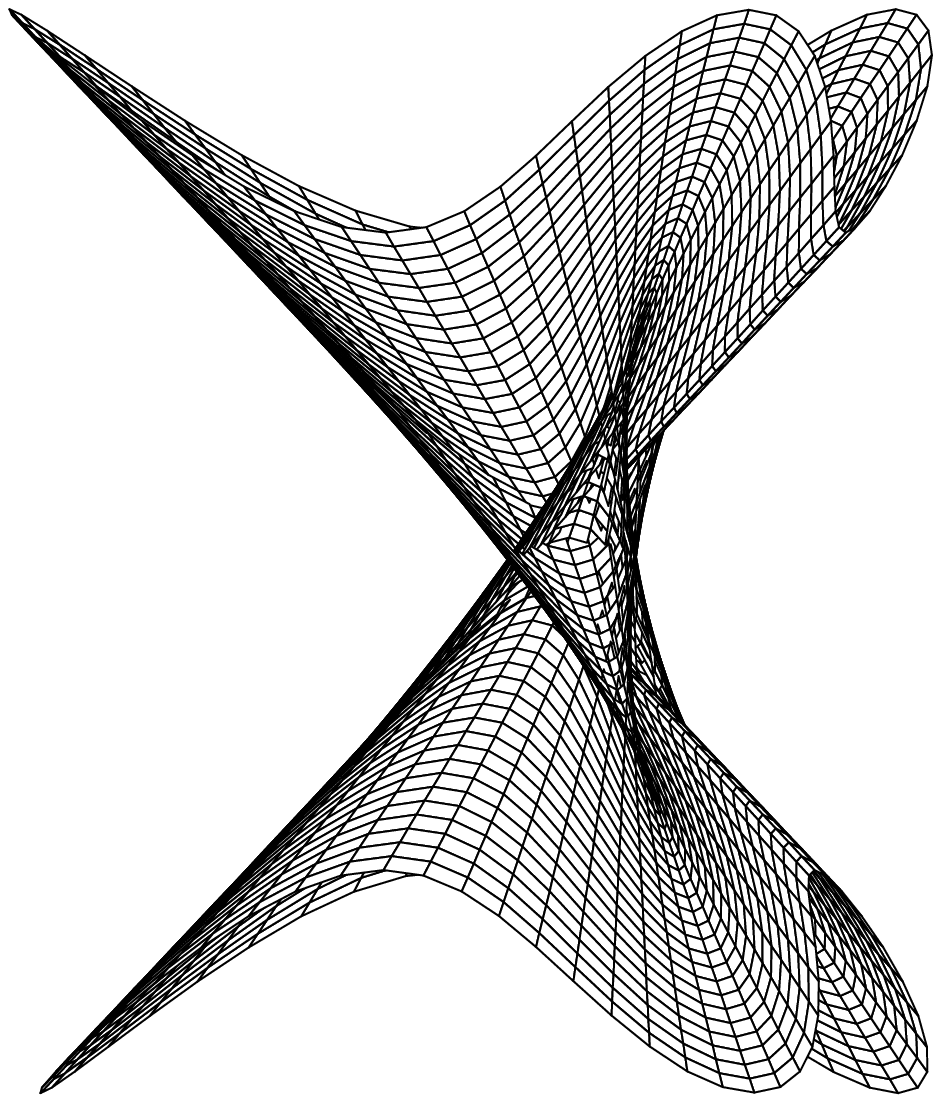}%
\\
Zoom on the zone C.
\label{chapt2f4}%
\end{center}

\newpage

\subsection{Computing the exponential mapping \label{appA1}}

Geodesics are trajectories of the Hamiltonian vector field $H$ associated to
the Hamiltonian $\mathcal{H}(\psi)$ on $T^{*}\mathcal{M}$ ($\pi:T^{*}%
\mathcal{M}\rightarrow\mathcal{M}$)
\[
\mathcal{H}(\psi)=\frac12\left(  \psi(F)^{2}+\psi(G)^{2}\right)  .
\]
The metric $(F,\,G)$ is in normal form (coordinates in $\mathcal{M}%
:\xi=(x,\,y,\,w)=(z,\,w)$)
\[%
\begin{array}
[c]{l}%
F=(1+y^{2}\,\beta)\frac\partial{\partial x}-x\,y\,\beta\,\frac\partial
{\partial y}+\frac y2(\gamma)\,\frac\partial{\partial w},\\
\\
G=(1+x^{2}\,\beta)\frac\partial{\partial y}-x\,y\,\beta\,\frac\partial
{\partial x}-\frac x2(\gamma)\,\frac\partial{\partial w},\\
\\
\beta(0,0)=b_{0}\\
\\
\gamma=(1+(x^{2}+y^{2})\beta)\,\smallint_{0}^{1}\frac{2\,t\,dt}{1+t^{2}%
(x^{2}+y^{2})\,\beta(t\,x,t\,y)}%
\end{array}
\]
Coordinates in the cotangent bundle are ($\widetilde{p},\,\widetilde{q}$,
$r)=(\widetilde{\zeta},\,r).$ $(x,\,y,\,w,\,\widetilde{p},\,\widetilde{q},r)$
have weight $1,\,1,\,2,-1,-1,-2$ respectively.
\[
\mathcal{H}^{-1}(\frac12)\cap\pi^{-1}(0)=\{\widetilde{p}^{2}(0)+\widetilde
{q}^{2}(0)=1\}.
\]
We set $\widetilde{p}(0)=\cos\varphi,\;\widetilde{q}(0)=\sin\varphi$. For
$r(0)\neq0$, we set $\rho=\dfrac1{r(0)},\;p=\dfrac{\widetilde{p}}%
r,\;q=\dfrac{\widetilde{q}}r,\;\zeta=(p,\,q).$ $s$ denotes the arclength and
$t$ the new time : $dt=r(s)\,ds$. $p,\,q$ have weight 1.

\smallskip\ 

One has, for $r(0)\neq0$ and $s$ small :
\begin{equation}
\frac{dz}{ds}=r(s)\,\frac{\partial\mathcal{H}}{\partial\widetilde{\zeta}%
}|_{r=1},\;\;\text{and\ }\;\;\frac d{ds}(\frac{\widetilde{\zeta}%
}r)=-r(s)\,\frac{\partial\mathcal{H}}{\partial z}|_{r=1} \label{A12}%
\end{equation}
Or,
\begin{equation}%
\begin{array}
[c]{l}%
\frac{dz}{dt}=\frac{\partial\mathcal{H}}{\partial\widetilde{\zeta}}%
|_{r=1},\;\;\;\;\;\frac{d\zeta}{dt}=-\frac{\partial\mathcal{H}}{\partial
z}|_{r=1}%
\end{array}
\label{A13}%
\end{equation}
For all $k$ (\ref{A13}) can be rewritten :
\begin{equation}
\frac{d(z,\,\zeta)}{dt}=A(z,\,\zeta)+\sum_{i=3}^{k}F_{i}(z,\,\zeta
,\,w)+O^{k+1}(z,\,\zeta,\,w), \label{A14}%
\end{equation}
where $F_{i}$ is homogeneous of degree $i,$ where $A$ is a linear operator
(corresponding to the Heisenberg sub-Riemannian metric), and where
$O^{k+1}(z,\,\zeta,\,w)$ has order $(k+1)$ with respect to the gradation :
\[
\deg x=\deg y=\deg p=\deg q=1,\;\deg w=2.
\]

\smallskip\ 

Also, $\dfrac{dw}{ds}=\dfrac{\partial\mathcal{H}}{\partial r}=r(s)\,\dfrac
{\partial\mathcal{H}}{\partial r}|_{r=1},$%
\begin{equation}
\dfrac{dw}{dt}=\dfrac{\partial\mathcal{H}}{\partial r}|_{r=1}. \label{A15}%
\end{equation}
This can be rewritten :
\begin{equation}
\dfrac{dw}{dt}=G_{2}(z,\,\zeta)+\sum_{i=4}^{k}G_{i}(z,\,\zeta,\,w)+O^{k+1}%
(z,\,\zeta,\,w), \label{A16}%
\end{equation}
where $G_{i}$ are homogeneous of degree $i,$ where $G_{2}$ corresponds to the
Heisenberg sub-Riemannian metric, and where $O^{k+1}(z,\,\zeta,\,w)$ has order
$(k+1)$ w.r.t. the gradation.

Initial conditions are
\begin{equation}
z(0)=0,\;w(0)=0,\;\zeta(0)=(\rho\,\cos\varphi,\,\rho\,\sin\varphi).
\label{A17}%
\end{equation}
Therefore,
\begin{equation}
\left\{
\begin{array}
[c]{l}%
(z,\,\zeta)=\rho\,(z_{1}(t,\,\varphi),\,\zeta_{1}(t,\,\varphi))+\sum
\limits_{i=3}^{k}\rho^{i}\,(z_{i}(t,\,\varphi),\,\zeta_{i}(t,\,\varphi
))+O(\rho^{k+1}),\\
\\
w=\rho^{2}\,w_{2}(t,\,\varphi)+\sum\limits_{i=4}^{k}\rho^{i}\,w_{i}%
(t,\,\varphi)+O(\rho^{k+1}).
\end{array}
\right.  \label{A17p}%
\end{equation}%
\begin{equation}
(z_{1},\,\zeta_{1})(t,\,\varphi)=e^{A\,t}\,(z(0),\,\zeta(0)),\;w_{2}%
(t,\varphi)=\int_{0}^{t}G_{2}(z_{1}(\tau,\varphi),\,\zeta_{1}(\tau
,\varphi))\,d\tau. \label{A18}%
\end{equation}
These last expressions can be easily computed. They give the exponential
mapping for the Heisenberg metric :
\begin{equation}
\left\{
\begin{array}
[c]{l}%
z_{1}(t,\,\varphi)=2\,\sin(\frac{t}{2})\,\left(  \cos(\varphi-\frac{t}%
{2}),\,\sin(\varphi-\frac{t}{2})\right)  ,\\
\\
\zeta_{1}(t,\,\varphi)=\cos(\frac{t}{2})\,\left(  \cos(\varphi-\frac{t}%
{2}),\,\sin(\varphi-\frac{t}{2})\right)  ,\\
\\
w_{2}(t)=\frac{1}{2}(t-\sin t).
\end{array}
\right.  \label{A19}%
\end{equation}
Also $F_{3}$ and $G_{4}$ don't depend on $w$, hence, setting $\Lambda
=(z,\,\zeta),$%
\begin{equation}
\left\{
\begin{array}
[c]{l}%
\Lambda_{3}(t,\varphi)=(z_{3},\zeta_{3})(t,\varphi)=\int_{0}^{t}%
e^{A\,(t-\tau)}\,F_{3}(\Lambda_{1}(\tau,\varphi))\,d\tau,\\
\\
w_{4}(t,\varphi)=\int_{0}^{t}\left(  \frac{\partial G_{2}}{\partial\Lambda
}(\Lambda_{1}(\tau,\varphi))\,.\,\Lambda_{3}(\tau,\,\varphi)+G_{4}(\Lambda
_{1}(\tau,\varphi))\right)  \,d\tau.
\end{array}
\right.  \label{A110}%
\end{equation}
The following terms are easily computed, on the same way. We give the
expressions that we shall need, $F_{4},\;G_{5}$ don't depend on $w$,
\begin{equation}
\left\{
\begin{array}
[c]{l}%
\Lambda_{4}(t,\varphi)=(z_{4},\zeta_{4})(t,\varphi)=\int_{0}^{t}%
e^{A\,(t-\tau)}\,F_{4}(\Lambda_{1}(\tau,\varphi))\,d\tau,\\
\\
w_{5}(t,\varphi)=\int_{0}^{t}\left(  \frac{\partial G_{2}}{\partial\Lambda
}(\Lambda_{1}(\tau,\varphi))\,\,\Lambda_{4}(\tau,\,\varphi)+G_{5}(\Lambda
_{1}(\tau,\varphi))\right)  \,d\tau.
\end{array}
\right.  \label{A111}%
\end{equation}%
\begin{equation}
\left\{
\begin{array}
[c]{l}%
\Lambda_{5}(t,\varphi)=\int_{0}^{t}e^{A\,(t-\tau)}\,(F_{5}(\Lambda_{1}%
(\tau,\varphi),\,w_{2}(\tau))+\frac{\partial F_{3}}{\partial\Lambda}%
(\Lambda_{1}(\tau,\varphi))\,\,\Lambda_{3}(\tau,\,\varphi))\,d\tau,\\
\\
w_{6}(t,\varphi)=\int_{0}^{t}\frac{1}{2}\frac{\partial^{2}G_{2}}%
{\partial\Lambda^{2}}\,.\,(\Lambda_{3}(\tau,\,\varphi),\Lambda_{3}%
(\tau,\,\varphi))+G_{6}(\Lambda_{1}(\tau,\varphi),w_{2}(\tau))+\\
\\
\;\;\;\;\;\;\;\;\;\;\;\frac{\partial G_{2}}{\partial\Lambda}(\Lambda_{1}%
(\tau,\varphi))\,.\,\Lambda_{5}(\tau,\,\varphi)\,+\frac{\partial G_{4}%
}{\partial\Lambda}(\Lambda_{1}(\tau,\varphi))\,.\,\Lambda_{3}(\tau
,\,\varphi))d\tau.
\end{array}
\right.  \label{A112}%
\end{equation}%
\begin{equation}
\left\{
\begin{array}
[c]{l}%
\Lambda_{6}(t,\varphi)=\int_{0}^{t}e^{A\,(t-\tau)}(F_{6}(\Lambda_{1}%
(\tau,\varphi),w_{2}(\tau))+\frac{\partial F_{3}}{\partial\Lambda}(\Lambda
_{1}(\tau,\varphi))\,\Lambda_{4}(\tau,\,\varphi)+\\
\\
\frac{\partial F_{4}}{\partial\Lambda}(\Lambda_{1}(\tau,\varphi))\,\Lambda
_{3}(\tau,\varphi))d\tau,\\
\\
\Lambda_{7}(t,\varphi)=\int_{0}^{t}e^{A\,(t-\tau)}\,(F_{7}(\Lambda_{1}%
(\tau,\varphi),\,w_{2}(\tau))+\frac{\partial F_{3}}{\partial\Lambda}%
(\Lambda_{1}(\tau,\varphi))\,\Lambda_{5}(\tau,\,\varphi)+\\
\\
\frac{\partial F_{4}}{\partial\Lambda}(\Lambda_{1}(\tau,\varphi))\,\Lambda
_{4}(\tau,\varphi)+\\
\\
\frac{1}{2}\frac{\partial^{2}F_{3}}{\partial\Lambda^{2}}(\Lambda_{1}%
(\tau,\varphi))\,(\Lambda_{3}(\tau,\varphi),\Lambda_{3}(\tau,\varphi
))+\frac{\partial F_{5}}{\partial\Lambda}(\Lambda_{1}(\tau,\varphi
))\,\Lambda_{3}(\tau,\varphi))d\tau,
\end{array}
\right.  \label{A113}%
\end{equation}
As we shall see, we will need to compute these values for $t=2\,\pi$ only.

\subsection{The exponential mapping in suspended form, the conjugate
loci\label{appA2}}

\smallskip\ %

\[
\mathcal{E}(t,\varphi):\left\{
\begin{array}
[c]{l}%
z(t,\rho,\,\varphi)=\rho\,z_{1}(t,\,\varphi)+\sum\limits_{i=3}^{7}\rho
^{i}\,z_{i}(t,\,\varphi)+O(\rho^{8}),\\
\\
w(t,\rho,\,\varphi)=\rho^{2}\,w_{2}(t)+\sum\limits_{i=4}^{6}\rho^{i}%
\,w_{i}(t,\,\varphi)+O(\rho^{7}).
\end{array}
\right.
\]
Integral expressions of $z_{i},\,w_{i}$ were computed in our appendix
\ref{appA1}.

\smallskip

Let us consider the variable $\epsilon=\pm1$, according to $w>0$ or $w<0$. Let
us set $t=2\pi\epsilon+\tau$ (the conjugate new time has an expansion
$t_{c}=2\pi\,\epsilon+O(\rho^{2})$, as is shown in our previous papers
\cite{CGK}, and this will appear again here in). Therefore, $\tau_{c}$, the
conjugate time $\tau$, has order $\rho^{2}.$%
\[
w_{2}(t)=\frac12(t-\sin t).\;\;w_{2}(\tau)=w_{2}^{\tau}(\tau)=\pi
\,\epsilon+\frac12(\tau-\sin\tau),
\]
hence, $w_{2}^{\tau}(0)=\pi\,\epsilon,\;w_{2}^{\tau^{\prime}}(0)=0,\;w_{2}%
^{\tau^{\prime\prime}}(0)=0.$ Therefore, we obtain the following important
fact :
\begin{equation}
w_{2}(\tau)=\pi\,\epsilon+O(\tau^{3})=w_{2}^{\tau}(0)+\tau^{3}\,\psi(\tau),
\label{A21}%
\end{equation}
for some smooth function $\psi.$

\smallskip\ 

Let us set :
\begin{equation}
h=\sqrt{\frac{\epsilon\,w}\pi}. \label{A22}%
\end{equation}%

\[%
\begin{tabular}
[c]{l}%
$w(\tau,\,\rho,\,\varphi)=\pi\,\epsilon\,\rho^{2}(1+O(\tau^{3})+\frac{1}%
{\pi\,\epsilon}(\rho^{2}\,w_{4}(\tau,\varphi)+\rho^{3}\,w_{5}(\tau
,\varphi)+\rho^{4}\,w_{6}(\tau,\varphi)+$\\
\\
$\rho^{5}\,w_{7}(\tau,\varphi)+\rho^{6}\,w_{8}(\tau,\varphi)+O(\rho^{7})).$%
\end{tabular}
\]
\[
h=\rho\,(1+O(\tau^{3})+\frac{1}{\pi\,\epsilon}(\rho^{2}\,w_{4}(\tau
,\varphi)+\rho^{3}\,w_{5}(\tau,\varphi)+\rho^{4}\,w_{6}(\tau,\varphi
)+O(\rho^{5})))^{\frac{1}{2}}.
\]
Using the implicit function theorem, we can solve this last equation in $\rho
$. After straightforward computations, we obtain :
\begin{equation}%
\begin{array}
[c]{l}%
\rho=h\,(1-h^{2}\frac{w_{4}^{\tau}(0)}{2\pi\epsilon}-h^{3}\frac{w_{5}^{\tau
}(0)}{2\pi\epsilon}+h^{4}(\frac{7w_{4}^{\tau}(0)^{2}}{8\pi^{2}}-\frac
{w_{6}^{\tau}(0)}{2\pi\epsilon})-w_{4}^{\tau^{\prime}}(0)\frac{h^{2}\tau}%
{2\pi\epsilon}+\\
\\
O(\tau^{3})+O(h^{5})+O(h^{2}\tau^{2})+O(h^{3}\tau)\,).
\end{array}
\label{A23}%
\end{equation}

\smallskip\ 

Let us set now
\begin{equation}
\tau=\theta\,h^{2}+\sigma\,h^{3}, \label{A24}%
\end{equation}
\ for some constant $\theta$.%

\begin{equation}%
\begin{tabular}
[c]{l}%
$\rho=h\,(1-h^{2}\epsilon\frac{w_{4}^{\tau}(0)}{2\pi}-h^{3}\epsilon\frac
{w_{5}^{\tau}(0)}{2\pi}+h^{4}(\frac{7w_{4}^{\tau}(0)^{2}}{8\pi^{2}}%
-\epsilon\frac{w_{6}^{\tau}(0)}{2\pi}-w_{4}^{\tau^{\prime}}(0)\frac
{\epsilon\,\theta}{2\pi})+$\\
\\
$c_{5}$ $h^{5}+c_{6}\,h^{6}+O(h^{7})).$\\
\\
where :\\
\\
$c_{5}=-\epsilon\,\frac{w_{7}^{\tau}(0)}{2\pi}+$ terms depending on the
invariants $\beta_{i},$ $(i\leq2)$\\
\\
$c_{6}=-\epsilon\,\frac{w_{8}^{\tau}(0)}{2\pi}+$ terms depending on the
invariants $\beta_{i},$ $(i\leq3)$%
\end{tabular}
\label{A25}%
\end{equation}
Otherwise,
\begin{equation}
z(\tau,\,\rho,\,\varphi)=\rho\,z_{1}(\tau,\,\varphi)+\sum_{i=3}^{9}\,\rho
^{i}\,z_{i}(\tau,\,\varphi)+O(\rho^{8}). \label{Aet}%
\end{equation}
We set
\[
z_{i}^{(k)}=\frac{d^{k}z_{i}}{d\tau^{k}},\;\;w_{i}^{(k)}=\frac{d^{k}w_{i}%
}{d\tau^{k}}.
\]
We also denote $z_{i}^{^{\prime}}$ for $z_{i}^{(1)}$ etc$\ldots$.

Replacing (\ref{A24}), (\ref{A25}), in this expression (\ref{Aet}), we obtain,
after tedious computations :
\begin{equation}
z(\sigma,h,\varphi)=\sum_{i=3}^{9}\,A_{i}\,h^{i}+O(h^{10}), \label{A26}%
\end{equation}%
\begin{equation}
\left\{
\begin{array}
[c]{l}%
A_{3}=(\theta\,\,z_{1}^{^{\prime}}+z_{3})|_{t=2\pi\,\epsilon},\;\;\;A_{4}%
=(\sigma\,z_{1}^{^{\prime}}+z_{4})|_{t=2\pi\,\epsilon},\\
\\
A_{5}=(\frac{-\epsilon}{2\pi}(\theta\,w_{4}\,z_{1}^{^{\prime}}+3\,w_{4}%
\,z_{3})+\frac{\theta^{2}}{2}z_{1}^{^{\prime\prime}}+\theta\,z_{3}^{^{\prime}%
}+z_{5})|_{t=2\pi\,\epsilon}\\
\\
A_{6}=\frac{1}{2\pi}(-\epsilon\,\sigma\,w_{4}z_{1}^{^{\prime}}-\epsilon
\,\theta w_{5}\,z_{1}^{^{\prime}}+2\pi\sigma\theta z_{1}^{^{\prime\prime}%
}-3\epsilon\,w_{5}z_{3}+2\,\pi\sigma z_{3}^{^{\prime}}-\\
\\
4\epsilon\,w_{4}\,z_{4}+2\,\pi\,\theta\,z_{4}^{^{\prime}}+2\,\pi
\,z_{6})|_{t=2\pi\,\epsilon}=\sigma\,A_{6}^{1}+A_{6}^{0},\\
\\
A_{7}=\sigma^{2}\,A_{7}^{2}+\sigma\,A_{7}^{1}+A_{7}^{0},\\
\\
A_{7}^{2}=(\frac{1}{2}z_{1}^{^{\prime\prime}})|_{t=2\,\pi}%
,\,\,\,\,\,\,\,\,\,A_{7}^{1}=\frac{1}{2\pi}(\epsilon\,w_{5}\,z_{1}^{^{\prime}%
}+2\,\pi\,z_{4}^{^{\prime}})|_{t=2\,\pi}\\
\\
A_{7}^{0}=(\frac{1}{2\,\pi}(\frac{7\theta}{4\,\pi}w_{4}^{2}-\epsilon
\,\theta^{2}\,w_{4}^{^{\prime}}-\epsilon\,\theta\,w_{6})\,z_{1}^{^{\prime}%
}-\frac{\epsilon}{2}\theta^{2}w_{4}\,z_{1}^{^{\prime\prime}})+3\,\pi
\,\theta^{3}\,z_{1}^{^{\prime\prime\prime}}+\\
\\
\frac{3}{2\,\,\pi}w_{4}^{2}\,z_{3}+3(\frac{7}{4\,\pi}w_{4}^{2}+\epsilon
\,\theta\,w_{4}^{^{\prime}}-\epsilon\,w_{6})\,z_{3}+3\,\epsilon\,\theta
\,w_{4}\,z_{3}+\pi\,\theta^{2}\,z_{3}^{^{\prime\prime}}-\\
\\
4\epsilon\,w_{5}\,z_{4}+5\epsilon\,w_{4}\,w_{5}+2\pi\,\theta\,z_{5}^{^{\prime
}}+2\pi\,z_{7})|_{t=2\pi}.\\
\\
A_{8}=\sigma\,A_{8}^{1}+A_{8}^{0},\,\,\,\,A_{9}=\sigma^{2}\,A_{9}^{2}%
+\sigma\,A_{9}^{1}+A_{9}^{0}.
\end{array}
\right.  \label{A27}%
\end{equation}

\smallskip\ \newline As we know from \cite{CGK}, the conjugate time $\tau_{c}$
is obtained for
\begin{equation}
\theta=-6\,\pi\,\epsilon\,\,b_{0} \label{A28}%
\end{equation}
To see this it is sufficient to compute $\theta$ for
\[%
\begin{array}
[c]{c}%
0=z_{1}^{\tau^{\prime}}(0,\varphi)\wedge(\theta\frac{\partial z_{1}%
^{\tau^{\prime}}(0,\varphi)}{\partial\varphi}+\frac{\partial z_{3}^{\tau
}(0,\varphi)}{\partial\varphi})\;\;\;\;\;\;\\
\\
\;\;\;=z_{1}^{t^{\prime}}(2\pi\epsilon,\varphi)\wedge(\theta\frac{\partial
z_{1}^{t^{\prime}}(2\pi\epsilon,\varphi)}{\partial\varphi}+\frac{\partial
z_{3}^{t}(2\pi\epsilon,\varphi)}{\partial\varphi}),
\end{array}
\]
but, $z_{1}^{t^{\prime}}(2\pi\epsilon,\varphi)\wedge\dfrac{\partial
z_{1}^{t^{\prime}}(2\pi\epsilon,\varphi)}{\partial\varphi}=1,\;\;z_{1}%
^{t^{\prime}}(2\pi\epsilon,\varphi)\wedge\dfrac{\partial z_{3}^{t}%
(2\pi\epsilon,\varphi)}{\partial\varphi}=6\,\pi\,\epsilon\,b_{0}.$

Also, we have :
\[
-6\,\pi\,\epsilon\,b_{0}\;z_{1}^{t^{\prime}}(2\pi\epsilon,\varphi)+z_{3}%
^{t}(2\pi\epsilon,\varphi)=0,
\]
as it is easily checked.

\medskip\ 

Therefore, we have the expression of the exponential mapping, in suspended
form, in a certain neighborhood $U$ of the conjugate locus at the source, for
$h$ (or $\rho$) small enough
\begin{equation}
z(\sigma,\,h,\,\varphi)=\sum_{i=4}^{9}\,A_{i}\,h^{i}+O(h^{10}), \label{A29}%
\end{equation}
where $A_{i}$ are given in (\ref{A27}).

We have to compute the expression of the conjugate time $\sigma_{c}$ in terms
of $h$ and $\varphi$. For this, we have to solve the following equation in
$\sigma$ :
\begin{equation}
\frac{\partial z}{\partial\sigma}(\sigma,\,h,\,\varphi)\wedge\frac{\partial
z}{\partial\varphi}(\sigma,\,h,\,\varphi)=0 \label{A210}%
\end{equation}
Solving this equation in $\sigma$, with the implicit function theorem gives :%

\[
\sigma_{c}=\underset{i=0}{\overset{5}{\sum}}h^{i}\,\sigma_{c}^{i}.
\]

\begin{remark}
$\sigma_{c}^{4}=-(z_{1}^{^{\prime}}(2\pi\epsilon)\wedge\frac{\partial A_{80}%
}{\partial\varphi})+$ terms depending on the invariant $\beta_{i},$ $(i\leq4)$
and $\sigma_{c}^{5}=-(z_{1}^{^{\prime}}(2\pi\epsilon)\wedge\frac{\partial
A_{90}}{\partial\varphi})+$ terms depending on the invariant $\beta_{i},$
$(i\leq5).$ We don't need the complete expression of those two coefficients.
\end{remark}

Now we give the expression of $\sigma_{c}^{i}$ $(i=1,2,3).$
\begin{equation}
\left\{
\begin{array}
[c]{l}%
\sigma_{c}^{0}=-(z_{1}^{^{\prime}}(2\pi\epsilon)\wedge\frac{\partial z_{4}%
}{\partial\varphi}),\,\,\,\,\,\,\,\,\sigma_{c}^{1}=-(z_{1}^{^{\prime}}%
(2\pi\epsilon)\wedge\frac{\partial A_{5}}{\partial\varphi})\\
\\
\sigma_{c}^{2}=(-(z_{1}^{^{\prime}}(2\pi\epsilon)\wedge\frac{\partial
A_{6}^{0}}{\partial\varphi})-(A_{6}^{1}\wedge\frac{\partial z_{4}}%
{\partial\varphi}(2\pi\epsilon))+\\
\\
(z_{1}^{^{\prime}}\wedge\frac{\partial z_{4}}{\partial\varphi})(2\pi
\epsilon)\,(\,A_{6}^{1}\wedge\frac{\partial z_{1}^{^{\prime}}}{\partial
\varphi}(2\pi\epsilon)+z_{1}^{^{\prime}}(2\pi\epsilon)\wedge\frac{\partial
A_{6}^{1}}{\partial\varphi}))\\
\\
\sigma_{c}^{3}=-(z_{1}^{^{\prime}}(2\pi\epsilon)\wedge\frac{\partial A_{7}%
^{0}}{\partial\varphi})-(A_{61}\wedge\frac{\partial A_{5}}{\partial\varphi
})-(A_{71}\wedge\frac{\partial z_{4}}{\partial\varphi}(2\pi\epsilon))-\\
\\
(z_{1}^{^{\prime}}\wedge\frac{\partial z_{4}}{\partial\varphi})^{2}%
(2\pi\epsilon)\,(2\,A_{7}^{2}\wedge\frac{\partial z_{1}^{^{\prime}}%
(2\pi\epsilon)}{\partial\varphi}+z_{1}^{^{\prime}}(2\pi\epsilon)\wedge
\frac{\partial A_{7}^{2}}{\partial\varphi})+\\
\\
(z_{1}^{^{\prime}}(2\pi\epsilon)\wedge\frac{\partial A_{5}}{\partial\varphi
})\,(A_{6}^{1}\wedge\frac{\partial z_{1}^{^{\prime}}(2\pi\epsilon)}%
{\partial\varphi}+z_{1}^{^{\prime}}(2\pi\epsilon)\wedge\frac{\partial
A_{6}^{1}}{\partial\varphi})+\\
\\
(z_{1}^{^{\prime}}\wedge\frac{\partial z_{4}}{\partial\varphi})(2\pi
\epsilon)(A_{7}^{1}\wedge\frac{\partial z_{1}^{^{\prime}}(2\pi\epsilon
)}{\partial\varphi}+z_{1}^{^{\prime}}(2\pi\epsilon)\wedge\frac{\partial
A_{7}^{1}}{\partial\varphi}+2\,A_{7}^{2}\wedge\frac{\partial z_{4}%
(2\pi\epsilon)}{\partial\varphi})
\end{array}
\right.  \label{A211}%
\end{equation}

It remains to replace the expression of $\sigma_{c}$ in the expression of $z $
in (\ref{A29}), to obtain the expansion of the conjugate locus :
\[
z_{c}(h,\varphi)=z(\sigma_{c},h,\varphi)=\sum_{i=4}^{9}\,f_{i}(\varphi
)\,h^{i}+O(h^{8}).
\]

The expressions of the $f_{i}$ $($for $i\leq4\leq6$ and $\epsilon=1)$ are
given in \ref{bf}. If we denote by $f_{i}^{-}$ the $f_{i}$ for $\epsilon=-1$,
we obtain :

$f_{i}^{-}=-f_{i}$ (for $4\leq i\leq6$)

The expressions of $f_{7}$ can be computed in the same way, just replacing
(\ref{A211}) in (\ref{A29}). This has been done with Mathematica. Let :
\[%
\begin{tabular}
[c]{l}%
$\beta_{4}=L_{44}\,(x^{2}+y^{2})^{2}+a_{44}\,(x^{4}+y^{4}-6\,x^{2}%
\,y^{2})+4\,b_{44}\,x\,y\,(x^{2}-y^{2})+$\\
$c_{44}\,(x^{4}-y^{4})-2\,d_{44}\,x\,y\,(x^{2}+y^{2}).$%
\end{tabular}
\]

Assuming $r_{2}=0,$ we obtain :%

\begin{equation}%
\begin{tabular}
[c]{l}%
1) $f_{7}=3\,\pi(-21\,c_{44}\,\cos(\varphi)+35\,a_{44}\,\cos(3\,\varphi
)-7\,c_{44}\,\cos(3\,\varphi)+21\,a_{44}\,\cos(5\,\varphi)+$\\
\\
$3\,r_{1}^{2}\,\cos(\varphi-2\,\theta_{1})+\,r_{1}^{2}\,\cos(3\,\varphi
-2\,\theta_{1})+21\,d_{44}\,\sin(\varphi)-35\,b_{44}\,\sin(3\,\varphi)+$\\
\\
$7\,d_{44}\,\sin(3\,\varphi)-21\,b_{44}\,\sin(5\,\varphi)-12\,\pi\,r_{1}%
^{2}\sin(\varphi-2\,\theta_{1})-12\,\pi\,r_{1}^{2}\sin(3\,\varphi
-2\,\theta_{1})),$\\
\\
$3\,\pi(21\,d_{44}\,\cos(\varphi)-35\,b_{44}\,\cos(3\,\varphi)-7\,d_{44}%
\,\cos(3\,\varphi)+21\,b_{44}\,\cos(5\,\varphi)-$\\
\\
$12\,\pi\,r_{1}^{2}\cos(\varphi-2\,\theta_{1})+12\,\pi\,r_{1}^{2}%
\cos(3\,\varphi-2\,\theta_{1})+21\,c_{44}\,\sin(\varphi)-35\,a_{44}%
\,\sin(3\,\varphi)-$\\
\\
$7\,c_{44}\,\sin(3\,\varphi)+21\,a_{44}\,\sin(5\,\varphi)-3\,r_{1}^{2}%
\,\sin(\varphi-2\,\theta_{1})+\,r_{1}^{2}\,\sin(3\,\varphi-2\,\theta_{1}))$\\
\\
2) $f_{7}+f_{7}^{-}=-144\,\pi^{2}\,\sin(2(\varphi-\theta_{1}))(\cos
(\varphi),\sin(\varphi))=\mathbf{d}_{7}.$%
\end{tabular}
\label{f7}%
\end{equation}

\smallskip\ 

It is not necessary to compute the expression of $f_{8}$, the only thing we
need is to know that
\begin{equation}
f_{8}(\varphi+\pi)=f_{8}(\varphi). \label{A213}%
\end{equation}
As we know, this is stated in \cite{Agr}, and we verify it by Mathematica.

\smallskip\ 

For $f_{9}$ the term $z_{9}$ appears in $A_{9}^{0}$ (formula (\ref{A27})), but
also, it comes through the term $A_{4}^{1}\,h^{4}$ of$\,\,z(\sigma
,\,h,\,\varphi)$ in (\ref{A27}) : $A_{4}^{1}=(\sigma\,z_{1}^{^{\prime}%
})\,(2\,\pi\,\epsilon).$

Therefore,the expression of $f_{9}$ contains also the term
\[
T=\sigma_{c}^{5}\,z_{1}^{^{\prime}}(2\,\pi\,\epsilon)\,h^{9}.
\]
This term will not play any role : the only thing that we need is to prove the
remark \ref{RC}.
\begin{equation}
\dfrac{\partial f_{6}}{\partial\varphi}=-180\,\left|  r_{3}\right|
\,\pi\,\sin(3\,\varphi+\theta_{3})\,(\cos\varphi,\,\sin\varphi),\;\;z_{1}%
^{^{\prime}}(2\pi\epsilon)=(\cos\varphi,\,\sin\varphi). \label{df6}%
\end{equation}
Hence, $\dfrac{\partial f_{6}}{\partial\varphi}\wedge f_{9}(\varphi
)=\dfrac{\partial f_{6}}{\partial\varphi}\wedge z_{9}(2\,\pi\,\epsilon)+$
other terms that depend on the invariants $\beta_{i},$ $(i\leq5).$

\smallskip\ 

We have :
\begin{equation}%
\begin{tabular}
[c]{l}%
$\dfrac{\partial f_{6}}{\partial\varphi}\wedge z_{7}(2\,\pi\,\epsilon
)=\dfrac{\partial f_{6}}{\partial\varphi}\wedge\int_{0}^{2\pi\epsilon
}\,e^{A(t-\tau)}\,F_{9}(\Lambda_{1}(\tau,\,\varphi),\,w_{2}(\tau))\,d\tau+$\\
\\
$\text{terms not depending on }\beta_{6}$%
\end{tabular}
\label{A214}%
\end{equation}

\smallskip\ 

This last expression has been computed by Mathematica, and we don't give it here.

We will give now the expression of the polynomial $P(\varphi)$ in the section
\ref{bf}.

Taking in account \ref{df6} and \ref{f7}, we see that $\dfrac{\partial f_{6}%
}{\partial\varphi}\wedge\mathbf{d}_{7}=0$. Therefor
\[
\dfrac{\partial f_{6}}{\partial\varphi}\wedge f_{7}(\varphi)=-\dfrac{\partial
f_{6}}{\partial\varphi}\wedge f_{7}^{-}(\varphi).
\]
$\,$ Hence (up to sign) the polynomial $P(\varphi)$ does not depend on
$\epsilon.$ Using Mathematica again we obtain :%

\[%
\begin{tabular}
[c]{l}%
$P(\varphi)=-7\,d_{44}\,\cos(2\,\varphi)+7\,b_{44}\,\cos(4\,\varphi
)-7\,c_{44}\,\sin(2\,\varphi)+7\,a_{44}\,\sin(4\,\varphi)+$\\
\\
$r_{1}^{2}\sin(2(\varphi-\theta_{1})).$%
\end{tabular}
\]

\end{document}